\documentclass[11pt,a4paper]{amsart}
\usepackage{amssymb}
\usepackage{amsmath}
\usepackage{amscd}
\usepackage[francais]{babel}
\usepackage{latexsym}
\usepackage[T1]{fontenc}

\newcommand{\spinendi}[1]{#1^* \!\!\otimes\! #1}
\newcommand{\spinendii}[2]{#1^*\!\!\otimes\! #2+#2^*\!\!\otimes\! #1}

\newcommand{\sanstr}[1]{ \{ #1 \}_0}
\newcommand{\pherm}[2]{ \langle #1,#2 \rangle }

\newcommand{\poids}{e^{2\delta t}}
\newcommand{\euc}{\mathrm{euc}}

\newcommand{\dxi}{\frac \partial {\partial {x_1}}}
\newcommand{\dyi}{\frac \partial {\partial {y_1}}}
\newcommand{\dxii}{\frac \partial {\partial {x_2}}}
\newcommand{\dyii}{\frac \partial {\partial {y_2}}}
\newcommand{\dzi}{\frac \partial{\partial {z_1}}}
\newcommand{\dzii}{\frac \partial{\partial {z_2}}}
\newcommand{\sh}{\operatorname{sh}}
\newcommand{\ch}{\operatorname{ch}}
\renewcommand{\tanh}{\operatorname{th}}
\newcommand{\argth}{\operatorname{arg\,th}}
\newcommand{\CHH}{\mathbb{C}\mathcal{H}^2}

\newcommand{\R}{\mathbb{R}}

\newcommand{\trace}{{\mathrm{tr}}}
\newcommand{\gmod}{{\hat g}}
\newcommand{\ghyp}{{ g ^{\mathcal{H}}}}
\newcommand{\ohyp}{{ \omega ^{\mathcal{H}}}}
\newcommand{\shyp}{{ s_{\mathcal{H}}}}

\newcommand{\gFS}{g^\mathrm{FS}}
\newcommand{\oFS}{\omega^\mathrm{FS}}

\newcommand{\delbar}{{\overline\partial}}
\newcommand{\delcp}{\partial^{h}}

\newcommand{\dt}{{\partial_t}}
\newcommand{\dtheta}{{\partial_\theta}}

\newcommand{\drho}{{\partial_\rho}}

\newcommand{\C}{\mathbb{C}}
\newcommand{\PP}{\mathbb{P}} 
\newcommand{\CP}{\mathbb{CP}}

\newcommand{\vol}{\mathrm{vol}}

\newcommand{\ffond}{\mathbb{I}}

\newcommand{\XX}{\mathcal{X}}

\newcommand{\jauge}{{\mathcal G}}
\newcommand{\Map}{{\mathrm Map}}

\newcommand{\Ric}{\mathrm{Ric}}
\newcommand{\Dirac}{\mathrm{D}}

\newcommand{\spinc}{spin^c}

\newcommand{\remarque}{\emph{Remarque}~: }

\renewcommand{\phi}{\varphi} \renewcommand{\Im}{{\mathrm{Im}}}
\renewcommand{\epsilon}{\varepsilon}

\newtheorem{prop}{Proposition}
\newtheorem{cor}[prop]{Corollaire}
\newtheorem{theo}[prop]{Théorème}
\newtheorem{resultat}{Résultat}

\newtheorem{lemme}[prop]{Lemme}
\newenvironment{demo}{\emph{Démonstration}. --- }{\hfill $\Box$ \medskip }
\newenvironment{demode}[1]{\emph{{Démonstration} #1: --- }}{\hfill $\Box$ \medskip}

\date{Décembre 2001}
\title{Rigidité d'Einstein du plan hyperbolique complexe}
\author{Yann Rollin}
\address{Yann Rollin \\ Department of Mathematics and Statistics \\
University of Edinburgh \\
James Clerk Maxwell Building \\
Kings Buildings \\
Mayfield Road \\
Edinburgh EH9 3JZ\\Scotland. }
\email{rollin@maths.ed.ac.uk}
\thanks{The author was supported by an EDGE grant, Research Training
Network HPRN-CT-2000-00101, European Human Potential Programme.}
\begin{document}
{\Huge \sc \bf\maketitle}

\begin{center}
\small
\begin{tabular}{p{13cm}}
{\sc Abstract.}\hspace{2mm}  We prove that every  Einstein metric on
$B^4\subset \C^2$ asymptotic to 
the Bergman metric is equal to it up to a diffeomorphism. We need a solution  of
Seiberg--Witten equations in this infinite 
volume setting. Therefore, and more generally, if $M^4$ is a manifold with a CR
boundary at infinity,  an adapted $\spinc$-structure which has a
non zero Kronheimer--Mrowka invariant and an asymptotically complex
hyperbolic Einstein metric, we produce a solution of  
Seiberg--Witten equations with an strong exponential decay property.
\end{tabular}
\end{center}
\bigskip

\section{Introduction}        
\bigskip

Les quotients compacts lisses de l'espace hyperbolique réel
$\R\mathcal H^4$ et du plan hyperbolique complexe $\CHH$ possèdent une unique
métrique d'Einstein à un difféomorphisme et une constante multiplicative
près. Dans le cas réel, ce résultat topologique est obtenu grâce à
des techniques d'entropie~\cite{BCG}. Dans le cas complexe, Le\,Brun
démontre ce résultat  par une voie très différente utilisant les équations de
Seiberg--Witten~\cite{L}.

En volume infini,
on sait d'après~\cite{CY} que les domaines strictement pseudo-convexes
de $\C^n$ possèdent une unique 
métrique complète de Kähler-Einstein à constante et biholomorphisme
près ; de façon analogue au cas compact, la question
difficile  de l'unicité de ces métriques  en 
tant que  métriques d'\emph{Einstein} se pose alors.    

Dans cet article, nous traitons le cas particulier de la boule unité  
$B^4$  
de $\C^2$ munie de la métrique de Bergmann $\ghyp$, qui est un modèle du  plan
hyperbolique complexe $\CHH$. On fixe un point de
la boule  et on note $t$ la distance à ce point relativement à la métrique $\ghyp$ ; nous démontrons le
théorème de rigidité suivant annoncé dans~\cite{R2}~:
\begin{theo}
\label{theo1}
Soit $g$ une métrique d'Einstein sur la boule unité $B^4\subset \C^2$ 
telle que près de l'infini on ait
\begin{equation}
\label{decasympt}
g= \ghyp + e^{-\delta t} h , \quad \mbox { $\delta
>0$}
\end{equation}
 et $h$ borné en norme $C^\infty$ relativement à $\ghyp$. Alors il existe un
difféomorphisme $f$ de $B^4$ tel que $f^*g = \ghyp$.
\end{theo}
\remarque en fait, l'hypothèse sur le comportement asymptotique de $g$
consiste, comme nous le verrons, à
ne considérer que  des métriques asymptotiquement hyperboliques
complexes dont l'infini conforme est celui de $\CHH$.

\medskip
La métrique hyperbolique complexe s'exprime classiquement en
coordonnées polaires comme suit~: notons $\rho$ la distance euclidienne
au centre de la boule $B^4$ d'où un difféomorphisme $B^4\setminus 0 \simeq ]0,1[\times S^3$. On obtient alors une
$1$-forme de contact $\eta$ sur $S^3$ définie par $\eta =
-Jd\rho /(2\rho)$ compatible avec le $S^1$-fibré de Hopf
$\pi:S^3\rightarrow \CP^1$. Alors en posant $t=\argth\rho$, on a 
$$\ghyp = dt^2+\sh^2(2t)\eta^2 + \sh^2(t)\gamma,
$$
où  $\gamma$ est la métrique de Carnot--Carathéodory sur la
distribution de contact  définie
par $\gamma = d\eta(\cdot, J\cdot) = \pi^*\gFS$ avec  $\gFS$ la métrique de
Fubini--Study sur  $\CP^1$. En fait $\gamma$ ne dépend du choix du
centre de la boule  qu'à un facteur conforme près et la classe
conforme $[\gamma]$ est appelée \emph{l'infini conforme}
de la métrique $g$.

 Pour chaque petite déformation de l'infini conforme
$[\pi^*\gFS]$, Biquard construit une métrique d'Einstein
asymptotiquement hyperbolique complexe sur la boule
\cite{B} ; l'hypothèse (\ref{decasympt}) du  théorème \ref{theo1}  est
donc bien entendu nécessaire. \medskip

La difficulté principale dans la démonstration du théorème~\ref{theo1} est d'exhiber une
solution des équations de Seiberg--Witten pour la métrique $g$, de
volume infini,  sans disposer d'un invariant défini a priori. Afin de
résoudre ce problème, nous approximons  dans la
section~\ref{collecone} la métrique $g$ par une suite
de métriques $g_\tau$ au bout desquelles on a adjoint un cône asymptotiquement
plat. 

La théorie de
Kronheimer et Mrowka \cite{KM} fournit alors un invariant de Seiberg--Witten non
nul relativement à la structure $\spinc$ canonique pour les métriques
$g_\tau$. On en déduit une suite de solutions $(A_\tau,\Phi_\tau)$ des
équations de Seiberg--Witten perturbées
\begin{eqnarray*}
\Dirac_{A_\tau}^{g_\tau} \Phi_\tau &=&0\\
F^+_{A_\tau\otimes B_\tau^{-1}} &= &q(\Phi_\tau) - q(\Phi_0) + \varpi_\tau;
\end{eqnarray*}
ici $B_\tau$ désigne une suite de connexions asymptotiquement plates tendant vers la connexion de
Chern $B$ induite par $\ghyp $ sur le fibré anti-canonique et
$(B,\Phi_0)$ est la solution standard des équations
pour la métrique hyperbolique complexe~; $\varpi_\tau$ est une suite de
$1$-formes autoduales à supports compacts qui convergent vers
$\varpi$ telle que  $F_B^+ = q(\Phi_0)+\varpi$. À la limite, la perturbation des 
équations se <<détache>> donc. Par ailleurs nous verrons que l'hypothèse sur le comportement asymptotique de $g$ implique  que
$\varpi$ possède une forte décroissance exponentielle.
Nous démontrerons alors, que quitte à faire des changements de jauge, et à
extraire une sous 
suite, nous pouvons faire converger $(\Phi_\tau,A_\tau)$ sur tout compact
vers une solution $(A,\Phi)$ des équations de Seiberg--Witten non
perturbées pour la métrique $g$  
\begin{eqnarray}
\label{SWnp}
\Dirac_{A}^{g} \Phi &=&0\\
F^+_A &= &q(\Phi),\nonumber
\end{eqnarray}
avec $a= A\otimes B^{-1}$ et $\phi= \Phi-\Phi_0$ dans $L^2_1$. Ce premier
résultat s'obtient via un contrôle uniforme sur l'\emph{énergie} des
solutions $(A_\tau,\Phi_\tau)$, section~\ref{secconvcompact}. 
Une méthode d'extraction similaire a  permis d'obtenir un résultat de
rigidité \cite{B2} pour les quotients de volume fini de
$\CHH$, mais il n'y a pas besoin de fixer d'infini conforme dans ce
cas et le résultat d'unicité concerne toutes les métriques
d'Einstein. En outre il est plus facile de faire converger les
solutions en volume fini, car la borne $C^0$ a priori obtenue sur la
partie spineur des solutions donne automatiquement un premier contrôle
$L^2$ et il n'y a donc pas besoin d'avoir recours à une énergie (cf. également \cite{R} pour une autre application en volume fini).

Le comportement  $L^2_1$   à l'infini est a priori insuffisant pour donner  à $F_A$ la  signification  cohomologique dont nous avons besoin (cf. section\ref{secdemo}).
Comme la métrique $g$ est d'Einstein, son
comportement asymptotique est alors précisé dans les travaux
d'Olivier Biquard et Marc Herzlich~\cite{BH}. 
Nous en déduirons,  en fixant une jauge de Coulomb près de l'infini, que $a$ et
$\phi$ en fait  des  $O(e^{-(4+\epsilon)t})$ avec $\epsilon>0$ ; pour
obtenir ce résultat, nous devrons analyser précisément  le
comportement asymptotique des équations
linéarisées et des laplaciens mis en jeux dans la
section~\ref{secbootstrapping} entre des espaces de Sobolev à poids.

Puis, la théorie de
Chern-Weil pour la signature et la caractéristique d'Euler
de~\cite{BH} nous permet  d'en déduire facilement le 
théorème~\ref{theo1} via une démonstration analogue à celle du cas
compact. 

Plus généralement, nous démontrons l'existence une solution des
équations de Seiberg--Witten dans les cas suivants :
\begin{theo}
\label{theoexistence}
Soit $M^4$, une variété orientée avec un bord à l'infini $Y^3$ muni d'une
structure CR et d'une $1$-forme de contact $\eta$ compatible avec
l'orientation. Soit $g$ une métrique d'Einstein asymptotiquement hyperbolique
complexe sur $M$, c'est à dire vérifiant  dans une trivialisation près l'infini
$$g = dt^2 + \sh^2(2t)\eta^2 + \sh^2(t)\gamma + O(e^{-\delta t}), \mbox{ 
$\delta >0$}
$$
où $O(e^{-\delta t})$ porte
également sur les dérivées successives et  $\gamma =
d\eta(\cdot,J\cdot)$.

Soit $\mathfrak s$ une
structure $\spinc$ sur $M$ compatible avec la structure de contact,
dont l'invariant de Kronheimer--Mrowka est non nul. Alors, les
équations de Seiberg--Witten correspondantes admettent 
une solution $(A,\Phi)$ telle que $A\otimes B^{-1}$, $\Phi
-\Phi_0$ et leurs dérivées sont des $O(e^{-(4+\epsilon)t})$ avec
$\epsilon >0$ suffisamment petit.
\end{theo}
 Dans l'énoncé, $(B,\Phi_0)$ est la solution standard des
équations de Seiberg--Witten pour une métrique formelle de Kähler--Einstein
asymptotiquement hyperbolique complexe $\bar g$ définie près de
l'infini et telle que $g-\bar g= O(e^{-4t})$ (cf. section \ref{secbootstrapping}).

 L'hypothèse que $g$ soit d'Einstein ne sert qu'à obtenir  des
renseignements supplémentaires sur sont comportement asymptotique ; on
peut par exemple supposer à la place dans le théorème
\ref{theoexistence}  que $g$ diffère  de $\bar g$ 
 par un $O(e^{-(4+\epsilon') t})$ avec $\epsilon'>0$.

 L'existence de cette solution devrait être justifiée plus
naturellement par une théorie de Seiberg--Witten développée
directement sur les variétés d'Einstein asymptotiquement
hyperboliques complexes.

Par ailleurs, le théorème \ref{theoexistence} permet de montrer des rapports
intéressants entre  la géométrie CR de $Y$  et l'existence de remplissage par des
métriques d'Einstein asymptotiquement hyperbolique complexe. En
 combinant  le théorème \ref{theoexistence} avec  \cite{BH}, Biquard obtient le résultat
suivant (cf. section \ref{ar}) : sous les hypothèses du théorème \ref{theoexistence}
, on a l'inégalité de Miyaoka--Yau
\begin{equation}
\label{MY}
0\leq \chi(M)-3\tau(M) +\nu(\partial_\infty M)  = \frac 1{8\pi^2}\int_M \left (
3|W_g^-|^2- |W_g^+|^2 + \frac {s_g^2}{24}\right
)\vol^g,
\end{equation}
où $W^\pm_g$ sont les composantes de la courbure de Weyl et  $\nu$ est l'invariant de la structure CR à l'infini de $M$ introduit dans  
\cite{BH} (où l'identité du membre de droite de (\ref{MY}) est démontrée). De plus on a égalité si et seulement si $g$  
est la métrique hyperbolique complexe. 

\bigskip
\begin{center}\bf{Remerciements}\end{center}

Je remercie chaleureusement Olivier Biquard pour les nombreuses
discussions que j'ai eues avec lui sur ce sujet et sur les métriques
d'Einstein asymptotiquement hyperboliques complexes. \bigskip

\section{Métrique hyperbolique complexe -- collage de cônes}
\label{collecone}
\bigskip

\subsection{Le plan hyperbolique complexe}
On munit $\C^3$ de la forme hermitienne $\langle z,z\rangle=  -\bar z_0z_0+  \bar z_1z_1 +\bar z_2z_2$ et on
définit classiquement le plan hyperbolique complexe par $\CHH =
\PP \{z\in \C^3, \langle z,z\rangle <0\} $.  
L'espace tangent à $\CHH$ en un point $z\in \C^3$ est défini par 
$T_z
\CHH= z^\perp=\{y\in \C^3, \langle z,y\rangle =0\}.$
 Plus
généralement on veut pouvoir considérer tout vecteur de $T\C^3$
comme un vecteur tangent à $\CHH$~: pour cela, il suffit de
considérer la projection orthogonale de $y\in T_z\C^3$ sur $T_z\CHH$. On
définit alors la métrique hyperbolique complexe sur $\CHH$ par
\begin{equation}
g_z ^{\mathcal H}(y,y) = \frac {\langle z,z\rangle \langle y,y\rangle -\langle
z,y\rangle \langle
y,z \rangle}{-\langle z,z\rangle ^2}.
\end{equation}
Elle est de Kähler-Einstein, sa courbure sectionnelle est comprise
entre $-1$ et $-4$ et sa forme de Kähler découle d'un potentiel
\begin{equation}
\label{potentiel}
\ohyp= - \frac i2 \partial\delbar \ln (- \langle z,z\rangle).
\end{equation}
\medskip 

\subsection{Le modèle de la boule}
Si on fixe la coordonnée homogène $z_0=1$, $\CHH$ est alors réalisé comme la
boule unité 
$$B^4=\{(z_1,z_2)\in \C^2, \; \rho(z_1,z_2) = \sqrt {|z_1|^2
  + |z_2|^2}<1\}\subset \C^2$$
 munie de la métrique de Bergmann 
\begin{equation}
\label{metrberg}
 \ghyp = \frac {\mathrm {euc}}{1-\rho ^2} + \frac{\rho ^2} {(1-\rho ^2)^2}\left [ (d\rho) ^2 +
  (Jd\rho)^2\right ].
\end{equation}
On note  $\pi : \C^2\setminus\{0\} \rightarrow \CP^1$ la projection
canonique ; on a des coordonnées polaires données par
\begin{eqnarray*}
\C^2\setminus\{0\} &\stackrel T \longrightarrow & ]0,+\infty [\times S^3 \\
(z_1,z_2) &\longmapsto &  (\rho ,u ),
\end{eqnarray*}
où $ S^3 = \partial B^4\simeq \partial_\infty \CHH$ est la sphère unité de
$\C^2$ et 
$$ u= \frac
1{\sqrt{|z_1|^2+|z_z|^2 }}(z_1,z_2).$$
Dans cette nouvelle carte, l'action de  $r e^{i\theta}\in \C ^*$ sur
$C^2\setminus \{0\}$ est donnée par 
$$ re ^{i\theta}\cdot(\rho,u) =  (r\rho, e ^{i\theta}\cdot u ),$$
où  $e ^{i\theta}\cdot u$ est l'action de Hopf de $S^1$ sur $S^3$~; on
note $\partial_r$ et  $\dtheta = i r\partial _r$ les champs de vecteurs sur $\C^*$
induits par les coordonnées polaires $re
^{i\theta}$. Ces champs de vecteurs sont envoyés via l'application tangente en $z=1$ à
l'action de $\C ^*$ respectivement sur $\rho\drho$ et $\dtheta= \rho J\drho$,
où par abus de langage $\dtheta$ désigne le champ de vecteurs
tangent à l'action de Hopf de
$S^1$ sur $S^3$.\medskip

On définit la métrique de Fubini-Study sur
$\CP^1$, à courbure sectionnelle $4$ et de volume $\pi$,
en coordonnées homogènes  par 
\begin{equation}
  \label{mFS}
  \gFS = \frac 1{(|z_1|^2+|z_2|^2) ^2} \left | z_2dz_1 -z_1dz_2 \right
  |^2
\end{equation}
et sa forme de Kähler est donnée par
égale à 
\begin{equation}
  \label{kFS}
   \oFS = \frac i2 \partial\delbar \ln( \rho^2 ).
\end{equation}
On vérifie aisément la formule suivante~:
\begin{equation}
  \label{ghyp2}
  \ghyp= \frac {|d\rho |^2 }{(1-\rho ^2)^2}  + \left (
\frac {2\rho}{1-\rho
    ^2} \right ) ^2 \left |\frac {Jd\rho}{2\rho}\right |^2 +\frac
{\rho ^2}{1-\rho ^2 } \pi ^* \gFS.
\end{equation}

Définissons la $1$-forme 
$\eta = -\frac {Jd\rho}{2\rho}.$
Il est facile de voir que cette $1$-forme est l'image réciproque d'une $1$-forme
sur $S^3$ dans les coordonnées polaires, qu'elle est invariante sous
l'action de Hopf et que 
$d\eta = \oFS$, où par abus de langage $\oFS$ désigne en fait $\pi^*\oFS$.
\medskip

\subsection{Raccordement d'un cône asymptotiquement plat.}
\label{racone}
Plus généralement on s'intéresse aux métriques  définies  à l'aide
d'un couple de fonctions  strictement positives $f= (f_1,f_2)$ ne dépendant
que de $\rho$,  et données par la formule 
$$g_f=  f_1(\rho)\euc + f_2(\rho) \left [|d\rho|^2+ |Jd\rho |^2\right ].
$$
On vérifie facilement que $g_f$ est compatible avec la structure complexe. Sa
forme de Kähler est donnée par
$$ \omega = \rho ^2 f_1 \oFS + (f_1 + f_2)d(\rho^2) \wedge \eta.
$$
On en déduit que $g_f$ est kählérienne relativement à
$J$ sous la condition 
$ \rho \drho f_1 = 2 f_2 $.
Alors,
en choisissant une fonction de la forme  $f_1(\rho)= C\rho^\alpha$ avec
$C,\alpha >0$, on vérifie  via  le  le changement de
variable $t ^2=C\rho ^{\alpha + 2}$ que  $\omega = d(t ^2 \eta)$~; on a
ainsi  défini une métrique de \emph{cône asymptotiquement
  plat} sur $\C^2\setminus \{0\}$.\medskip 

Soit la famille de  de fonctions $(f_s)_{s\in]0,1[ }$ définie par
\begin{eqnarray*}
  f_s (\rho) &= \frac 1{1-\rho ^2} & \mbox{si $\rho<s$,} \\
  f_s (\rho) &= C_s\rho ^{\alpha_s} & \mbox{sinon,} 
\end{eqnarray*}
avec des constantes $C_s,\alpha_s >0$ choisies de telle façon que la fonction
$f_s$ soit $C ^1$ par morceaux (on peut également choisir des fonctions $C
^\infty$ très proches de cette construction). On note
 $g_s :=g_f $
pour $f={(f_s,\frac 1 2 \rho \drho  f_s) }$.
On vient de définir une suite de métriques kählériennes asymptotiquement
plates $g_s$  sur $\C^2$, telle que
$$\lim_{s \rightarrow 1}  g_s =\ghyp $$
sur tout compact de $B^4$.
\medskip 

\subsection{Métriques asymptotiquement hyperboliques complexes}
\label{subhc}
En faisant le changement de variable
$t=\argth \rho \in ]0,+\infty[$, il vient 
\begin{equation}
  \label{metrhyp}
  \ghyp  = dt^2 + \sh^2(2t)\eta ^2 + \sh^2 (t) \pi ^* \gFS.
\end{equation}
De même, les métriques $g_s$ sont de la forme
\begin{equation}
  \label{metrgen}
  g^{\mathcal H}_\tau = dt^2 + \xi_1 ^2\eta ^2 + \xi_2^2  \pi ^* \gFS,
\end{equation}
où  $\xi_1$ et $\xi_2$ sont des fonctions de $t$ vérifiant
\begin{equation}
  \label{condkahl}
  \xi_1 =2\xi_2\dt \xi_2,
\end{equation}
et affines pour $t$ assez grand.

Plus généralement, soit $M^4$  une variété orientée munie d'un bord à l'infini
$Y^3$ doté d'une structure CR et d'une $1$-forme de contact $\eta$ compatible avec
l'orientation. On note $\gamma = d\eta(\cdot,J\cdot)$ la métrique de
Carnot-Carathéodory sur~$Y$.

On dit qu'une métrique  est
\emph{asymptotiquement hyperbolique complexe} si il existe une
trivialisation $[0,+\infty[ \times Y$ de $M$ près de l'infini telle
que 
$$g=\gmod + O(e^{-\delta t}) \quad \mbox{ avec $\delta >0$ \; et }\quad \gmod = dt^2 + \sh^2 (2t)\eta ^2 + \sh ^2 (t) \gamma .
$$
Précisons que $O(e^{-\delta t})$ concerne également toutes les dérivées.
De façon analogue au cas du plan hyperbolique, 
on définit une suite de fonctions $(\mu_\tau)_{\tau > 0}$ comme
suit~: pour $t \leq \tau$ on pose $\mu_\tau(t)= \sh t $ puis on   fait
un recollement $C^2$ par morceaux avec la fonction
affine $ (t-\tau)\ch\tau  + sh\tau $ qu'on approxime par une
fonction $C^ \infty$ ; nous pourrions nous contenter de la
construction singulière en donnant un sens aux dérivées de $\mu_\tau$ en tant
que distributions.  Pour régulariser $\mu_\tau$, il suffit de faire
décroître  brutalement 
$\partial_t ^2\mu_\tau$ de $\sh (\tau -\epsilon )$ à $0$ quitte
à rendre $\partial_ t
^3\mu_\tau $ très négatif entre $\tau-\epsilon $ et $\tau$.
On en déduit  une suite  de métriques
asymptotiquement plates $\gmod_\tau$ par la
formule~(\ref{metrgen}), où $\xi_2=\mu_\tau$ et $\xi_1 =
2\mu_\tau\dt\mu_\tau$.  Nous choisirons $\epsilon$ ensuite afin que $\mu_\tau$
soit arbitrairement proche du modèle singulier en norme $C^2$. 

Puis on étend la structure $CR$ en une structure presque complexe au
voisinage de l'infini en posant $J_\tau dt = -\xi_{1,\tau} \eta$.
On constate que la métrique hermitienne $(\gmod_\tau,J_\tau)$ est
presque kählérienne.

Pour étudier les métrique $\gmod_\tau$, on utilise le champ de Reeb
$R$ défini par $i_R\eta=1$ et $i_R d\eta = 0$ et on choisit une base locale
$Y_1,Y_2$ de champs de vecteurs sur $Y^3$ de la distribution de
contact $\ker\eta$ avec $Y_2=J Y_1$. On note $x= \xi_1^{-1}R$, $y_i =
\xi_2^{-1}Y_i$.  Il est facile de voir que
$$
[\dt,x] = - \frac{\dt \xi_1} {\xi_1} x  ,\quad
[\dt, y_i] =  - \frac{\dt \xi_2} {\xi_2} y_i  ,\quad
[x, y_i] = O({\xi_2^{-1}}) ,\quad
[y_1,y_2]  = O(\xi_2^{-1}) ,$$
où les termes en $O(\xi_2^{-1})$ portent sur également sur les dérivées. 
Il en résulte que dans cette trivialisation, le tenseur de Nijenhuis,
la connexion de Levi--Civita, la connexion de Chern de
$(\gmod_\tau,J_\tau)$ et leurs courbures sont égales à celle obtenues
en partant de la métrique hyperbolique complexe $\ghyp$ à un
$O(\xi_2^{-1})$ près. On pourra grâce à cette remarque se contenter la plupart du temps de faire des
calculs dans ce cas modèle  
afin d'en déduire des renseignements sur
 les métriques presque kählériennes $(\gmod_\tau,J_\tau)$ et $(\gmod,J)$.  \medskip

On construit maintenant une suite de métriques approximant une métrique
asymptotiquement hyperbolique $g$ comme suit :
soit $\chi(t)$ une fonction lisse valant $0$ pour  $t\leq 0$, $1$ pour
$t\geq \frac 12$ et strictement croissante entre $0$ et $\frac 12$ ;
on pose 
\begin{equation}
\label{chitau}
\chi_\tau(t)= \chi(t-\tau+1)
\end{equation}
 puis on définit la suite de métriques
\begin{equation}
\label{defgtau}
g_\tau = \chi_\tau \gmod_\tau + (1-\chi_\tau)g.
\end{equation}

\subsection{Bulle de courbure positive.}
\label{secbulle}
Si on se contentait d'un recollement $C^2$ par morceau pour construire la fonction
$\xi_{2,\tau}$, la métrique $g_\tau$ aurait une singularité dans sa
courbure (car on pose $\xi_1= \dt(\xi_2^2)$)~; après
avoir régularisé la fonction, ce phénomène
se manifeste toujours par une bulle de courbure positive :
 \begin{lemme}
\label{courbchern}
   La courbure $F_{B_\tau}$ de la connexion de Chern du fibré anti-canonique
   associée à la métrique $\gmod_\tau$ et la structure presque
   complexe $J_\tau$ se décompose en 
   \begin{equation}
  F_{B_\tau} = F^b_\tau  - i \chi_\tau \;\frac {\partial_t^2\xi_1}{\xi_1}\;dt \wedge \xi_1
\eta ,        
   \end{equation}
où $|F^b_\tau|$ est uniformément bornée. De plus, on a
   \begin{equation}
  2 i\Lambda F_{B\tau} \geq  \shyp + O(\xi_2^{-t}) \mbox{ pour $t\leq \tau$
   et } \quad
  2 i\Lambda F_{B\tau} =  O(\xi_2^{-1}) \mbox { pour $t\geq \tau$.}
   \end{equation}
 \end{lemme}\medskip
On en déduit immédiatement le corollaire suivant dans le cas de $\CHH$. Si
   $g$ est une métrique asymptotiquement hyperbolique complexe, la
   courbure des métriques $\gmod_\tau$  diffère de celle calculée dans
   le cas de $\CHH$ par  un $O(\xi_{2,\tau}^{-1})$ et le corollaire
   reste valable dans ce cas.
\begin{cor}
  La courbure scalaire des métriques $(g_\tau)$ se décompose
  en
\begin{equation}
s= s_\tau ^b  -2 \chi_\tau \frac {\partial ^2 \xi _{1} }{\xi_{1}},
\end{equation}
où $s_\tau ^b$ est bornée indépendamment de $\tau$.
\end{cor}

 \begin{demode}{du lemme \ref{courbchern}} La métrique et la structure
presque complexe du fibré anti-canonique sont explicites près de
l'infini ; la connexion de Chern est donc elle même explicite ainsi
que sa courbure et on en déduit le lemme. 

Comme nous l'avons remarqué, il suffit d'étudier  le cas où $\gmod$ est la métrique
hyperbolique complexe et on en déduit la courbure de la connexion de
Chern à un $O(\xi^{-1}_{2,\tau})$ près. Dans ce cas, 
les approximations $\gmod_\tau$ sont
kählériennes 
(cf. \ref{racone}), ce qui simplifie beaucoup les calculs.
On commence par choisir une section  holomorphe du fibré canonique
$\Lambda^{2,0} T\C^2$~:
$$\varpi = \dzi\wedge\dzii,
$$
avec $z_1=x_1+iy_1$, $z_2= x_2+iy_2$
$$\dzi =\frac 12 \left (\dxi -i \dyi \right ) \quad \mbox{ et
  } \quad \dzii =\frac 12  \left ( \dxii -i \dyii \right ).
$$
Posons
$$u= \bar z_2 \dzi -\bar z_1\dzii, \quad v=z_1\dzi+z_2\dzii.
$$
Puisque 
$$\gmod_\tau = \frac{f_2}{\rho^4}|\bar z_1 dz_1 +\bar z_2 dz_2|^2 + \frac{f_1}{\rho^2} | z_2 dz_1 - z_1 dz_2|^2,
$$
on en déduit que les vecteurs $u$ et $v$ sont orthogonaux, que
$$|u|^2 = \rho ^2 f_1= \xi_2 ^2,\quad |v|^2 = f_2 = \frac 14 \xi_1 ^2
$$
et finalement 
$$ |\varpi|=\frac 1{\rho ^2}|u||v| =\frac 1\rho \sqrt{f_1f_2} = \frac
{\xi_1\xi_2}{2\rho ^2(t)}.
$$
Alors
$$\frac 12 F_{B_\tau}=\partial\delbar\ln |\varpi |= \partial\delbar \ln \xi_{1}+ \partial\delbar \ln \xi_{2} -
\partial\delbar \ln \rho ^2,
$$
or $\partial\delbar \ln \rho ^2 = -2id\eta$. Par ailleurs, on calcule
\begin{align*}
 \delbar \ln \xi_1 & = (\dt \ln \xi_1) dt^{0,1}  =  (\dt \ln
\xi_1)\frac { dt+ iJdt}2  = \frac 12 (\dt \ln \xi_1) (dt -i\xi_1\eta).
\end{align*}
Puisque $\partial\delbar =d\delbar$, il vient
$$\partial\delbar \ln \xi_1 = -\frac i2 \frac {\partial_t^2\xi_1}{\xi_1}  dt\wedge
\xi_1\eta -\frac i2\left( \frac{\xi_1}{\xi_2} {\dt \ln \xi_1}\right ) \xi_2 d\eta
$$
d'où
$$\partial\delbar \ln \xi_2 = -\frac i2 \left ( \partial_t^2\ln \xi_2 + (\dt
  \ln \xi_1)( \dt \ln \xi_2 )\right )  dt\wedge
\xi_1\eta -\frac i2\left (\frac {\xi_1}{\xi_2}\dt \ln \xi_2 \right ) \xi_2 d\eta.
$$
En regroupant les termes, on en déduit la formule
\begin{align}
  \label{cric}
  F_{B_\tau} = - \left [\frac {\partial_t^2\xi_1}{\xi_1} +\partial_t^2\ln \xi_2 + (\dt
  \ln \xi_1)( \dt \ln \xi_2 ) \right ]\; & i dt\wedge
\xi_1\eta  \nonumber \\
-\left [ \frac {\xi_1}{\xi_2}\dt \ln (\xi_1\xi_2)\right ] \; & i\xi_2d\eta.
\end{align}
 On en déduit le lemme par construction des fonctions $\xi_1$ et $\xi_2$.
\end{demode}
\bigskip

\section{Équations de Seiberg--Witten --- Énergie}
\label{secconvcompact}
\bigskip

Supposons $M$ munie d'une structure $\spinc$  adaptée à la structure de
contact  notée $\mathfrak s$~; ceci signifie que près de l'infini, $\mathfrak s$ est isomorphe à la structure $\spinc$ standard $\mathfrak s_0$ induite par la structure symplectique $\omega^{\gmod_\tau}$. En particulier, le fibré de spineurs 
$W=W^+\oplus W^-$ est isomorphe près de l'infini à 
\begin{align}
\label{scano}
  W ^+ = & \Lambda ^{0,0}M \oplus \Lambda ^{0,2} M, \\
  W ^- = & \Lambda ^{0,1} M\nonumber
\end{align}
où $\Lambda^{p,q}M$
est le fibré des formes de type $(p,q)$ par rapport à 
$J_\tau$. Le \emph{fibré déterminant} de cette structure $\spinc$ est donné
près de l'infini par  
$$L = K^{-1}_{J_\tau} =\Lambda ^2 W^+ =\Lambda ^2 W^-.$$
 Comme les structures presque complexes $J_\tau$ sont
homotopes (et même égales pour $t$ suffisamment petit), les fibrés de $(p,q)$-formes associés sont isomorphes et on
omettra généralement de préciser la dépendance en $\tau$.
On supposera en outre que la classe de cohomologie près de l'infini
$[\omega^{\gmod_\tau}]$ se prolonge à l'intérieur de $M$. Dans le cas
où $M$ est le plan complexe cette condition est évidemment vérifiée
mais plus généralement, nous avons besoin de faire cette hypothèse
pour obtenir appliquer la théorie de Kronheimer--Mrowka
(cf. démonstration de la proposition \ref{propenergie}).

Soit $\Phi_0$ le champ un spineurs de $W^+$ défini près de l'infini
par  $\Phi_0=(\sqrt 
{-\shyp}, 0)$ et $B_\tau$ la connexion de Chern de $L$ 
associée à la métrique hermitienne déduite de $\gmod_\tau$.
On a la formule bien connue pour l'opérateur de Dirac défini par $\gmod_\tau$ et $B_\tau$
$$\Dirac^{\gmod_\tau}_{B_\tau }(\alpha,\gamma)= \sqrt 2 \left (\delbar
\alpha +\delbar ^*\gamma 
\right );
$$
en particulier $\Dirac^{\gmod_\tau}_{B_\tau } \Phi_0 = 0$. La suite de
connexions $B_\tau$ tend par construction sur tout compact de $M$ vers
la connexion de Chern de $\gmod$.\medskip

 Kronheimer et Mrowka
développent dans \cite{KM} une théorie de Seiberg-Witten pour des
métriques presque kählériennes asymptotiquement plates --
ce qui est le cas des métriques $g_\tau$ -- et des équations de
Seiberg-Witten perturbées près de l'infini, de la forme
\begin{eqnarray}
  \label{SWp1}
  \Dirac ^ {g_\tau} _{A} \Phi & =  & 0  \\
  \label{SWp2}
  F ^{+_{g_\tau}}_{A}- q(\Phi) & = &  F^{+_{g_\tau}}_{ B_\tau}-q(\Phi_0),
\end{eqnarray}
avec $\Phi$ un champ de spineurs de $W^+$,  $A$ une connexion
unitaire sur $L$ et $q(\Phi)=\{\Phi ^* \otimes\Phi\}_0$
l'endomorphisme hermitien sans trace identifié à la $2$-forme
autoduale imaginaire pure   $F ^{+_{g_\tau}}_{A\otimes B_\tau ^{-1}}$
via le produit de Clifford.

Kronheimer et Mrowka montrent que l'espace des modules des solutions 
 des équations $(A,\Phi)$
modulo l'action du groupe de jauge $\mathcal G= \Map(M,S^1)$,
 donnée par   
$$ u\cdot(A,\Phi)= (A +2\frac{du}u,u^{-1}\Phi),
$$
est fini ; en comptant algébriquement le nombre de points de l'espace
des modules 
obtient l'invariant de Kronheimer--Mrowka $SW(\mathfrak s)$, qui ne
dépend pas de la métrique $g_\tau$.

 En particulier, s'il existe une métrique de
Kähler-Einstein asymptotiquement symétrique, sa forme de Kähler
définit un remplissage symplectique de $M$. En notant  $\mathfrak s_0$ 
la structure $\spinc$ canonique associée, on à $SW(\mathfrak s_ 0)=1$
d'après le théorème  $1.1$ de \cite{KM}.

Lorsque l'invariant de Seiberg--Witten est non nul, on obtient donc  une suite
de solutions des équations de Seiberg-Witten $(A_\tau,\Phi_\tau)$,
pour chaque métrique $g_\tau$ telles que
\begin{equation}
  \label{sol}
  A_\tau = B_\tau + a_\tau \quad\mbox { avec } a_\tau \in L^2_l(g_\tau),\quad
  \Phi_\tau = \Phi_0 +\phi_\tau \quad\mbox { avec } \phi_\tau\in L^2_{l,B_\tau}(g_\tau), 
\end{equation}
avec $l$ arbitrairement grand. Par ellipticité,on peut supposer que
les solutions $(A_\tau,\Phi_\tau)$ sont lisses.  \medskip 

On a l'identité bien
connue pour les métriques  kählériennes à courbure scalaire
constante $F_{B}^+=q(\Phi_0)$. On s'attend à ce que cette équation
soit vérifiée également par les métriques d'Einstein
asymptotiquement hyperboliques complexes, à un terme correctif
près  à forte décroissance exponentielle. Nous verrons à la section 
\ref{secbootstrapping} comment modifier la famille de perturbations
afin qu'elle s'évanouisse pour les équations limites.

\subsection{Contrôle $C^0$}
En appliquant le principe du maximum on a une première estimée $C^0$  uniforme a
priori sur les spineurs $\Phi_\tau$.
\begin{lemme}
\label{controleC0}
Il existe une constante $c$ telle que pour toute solution  $(A,\Phi)$ des
  équations de Seiberg-Witten, relativement à une métrique
  $g_\tau$, on ait $|\Phi|<c$ en tout point de $M$.
\end{lemme}
\begin{demo}
Puisque $\Phi_\tau - \Phi_0\in L^2_l$, en choisissant $l\geq 3$, $|\Phi_\tau - \Phi_0|$
tend vers $0$ à l'infini sur le cône asymptotiquement plat d'après
l'inclusion $L^2_3\subset C^0$~;
si $|\Phi_\tau| > |\Phi_0|$ en un point, alors le maximum
de $|\Phi|$ est atteint. En un maximum, 
\begin{align*}
  0 \leq & \frac 12\Delta|\Phi|^2
  = \langle \nabla_{A_\tau} ^*
  \nabla_{A_\tau}\Phi_\tau,\Phi_\tau\rangle -\langle \nabla_{A_\tau}
  \Phi_\tau, 
  \nabla_{A_\tau}\Phi_\tau\rangle \\
\leq  & \langle \nabla_{A_\tau} ^*
  \nabla_{A_\tau}\Phi_\tau,\Phi_\tau\rangle =  \langle \Dirac_{A_\tau} ^2\Phi_\tau,\Phi_\tau\rangle - \frac 12
  \langle F_{A_\tau}^+ \cdot \Phi_\tau,\Phi_\tau\rangle -
  \frac s4 |\Phi_\tau|^2,
\end{align*}
et en utilisant les équations, on en déduit que 
$$0\leq - \frac 14 |\Phi_\tau| ^4 - \frac 12 \langle (F_{B_\tau}^+ - q(\Phi_0))
\cdot \Phi_\tau,\Phi_\tau \rangle - \frac s4 |\Phi_\tau|^2. 
$$
On définit l'opérateur linéaire uniformément borné $$P_\tau ^b= \frac {s_\tau ^b}4 +
\frac 12 (F_{\tau}^b) ^ + - \frac 12 q(\Phi_0)$$
et l'opérateur
$$P_\tau = -\chi_\tau \frac 12 \frac{\partial_t^2 \xi_{1,\tau }}{\xi_{1,\tau}} -
\chi_\tau \frac i2 \frac {\partial_t^2\xi_{1,\tau}}{\xi_{1,\tau}}(dt\wedge \xi_{1,\tau}\eta
) ^+ .$$
Puisque $(dt\wedge \xi_{1,\tau} \eta
) ^+$ agit par multiplication de Clifford avec valeurs propres $\pm i$,
on en déduit que $P_\tau$ agit avec pour valeurs propres $0$ et
$-\chi_\tau  \xi_{1,\tau} ^{-1}\partial_t ^2 \xi_{1,\tau}$  et que la
partie négative des valeurs propres de $P_\tau 
+P_\tau ^b$ est uniformément minorée. Par le principe du maximum,
il en résulte que $|\Phi_\tau|$ est bornée par une constante
indépendante de $\tau$.
\end{demo}

Nous précisons maintenant les conventions utilisées dans le cadre
presque complexe. 
La connexion de Chern induite par $\gmod_\tau$ sur $W^+=\Lambda^ {0,0}
\oplus \Lambda ^{0,2}$ est donnée par 
$\nabla = d\oplus \nabla
_{B}$, où  $B$ est la connexion unitaire induite par la métrique
sur $K^{-1}=\Lambda ^{0,2}$ et $\nabla _{B}$ est la connexion hermitienne 
associée à $B$.  

Munissons le fibré en droites complexes trivial $\C$ de la connexion 
$d_a = d + a$, où $a$ est une $1$-forme imaginaire pure~; on en
déduit une connexion hermitienne 
$d_{\hat a } = d + \hat a  = d + \frac 12 a$
sur $\C^{1/2}$ qui est la racine carrée de $d_a$.
Si on twiste le fibré des spineurs $W^+$ par  $\C^{1/2}$, on en
déduit une connexion 
$\nabla _{\hat a } = d_{\hat a}\oplus
(\nabla_{B}\otimes d_{\hat a})$
 sur $W^+\otimes \C^{1/2}=W^+$. On note alors
$\delbar _{\hat a} = \nabla ^{0,1}_{\hat a}$~;
la
connexion $\nabla _{\hat a}$ induit une connexion sur $L= \Lambda ^2
W^+ = \Lambda ^{0,0}\otimes \Lambda ^{0,2}$ égale à 
$$ \nabla _{B} \otimes d_{\hat a}\otimes d_{\hat a}  = \nabla _B
\otimes d_{a} = \nabla _{B+a}\;;
$$ 
on note alors $A=B+a$ la connexion unitaire correspondante.
L'opérateur de Dirac $\Dirac_A$ associé à
la connexion $A$ et la métrique hermitienne $\gmod_\tau$ est donc
donné par
$$ \Dirac_A^{\gmod_\tau} \Phi = \sqrt 2 \left (\delbar_{\hat a} \alpha +
  \delbar_{\hat a} ^*\gamma \right ) ,$$

et les équations de Seiberg-Witten s'écrivent (en notant $F_a := da$) 
\begin{eqnarray}
  \label{SWcmplx}
 \delbar_{\hat a}
\alpha +\delbar_{\hat a} ^*\gamma & = & 0\\
   \Lambda F_a  & = & \frac i2(|\alpha |^2 - |\gamma|^2 + \shyp)  \\
  F_a ^{0,2} & = & \frac {\overline\alpha \gamma}2 .
\end{eqnarray}

Rappelons  les formules de Weitzenböck  très utiles (cf. \cite{Ko}) 
\begin{align}
  \label{Weitz1}
\delbar ^*_{\hat a}\delbar _{\hat a} \alpha &=  \frac 12 (\nabla
  ^*_{\hat a}\nabla_{\hat a} \alpha - i \Lambda F_{\hat a}\alpha),  \\
  \label{Weitz2}
\delbar _{\hat a}\delbar^* _{\hat a} \gamma &=  \frac 12 (\nabla
  ^*_{\hat a}\nabla_{\hat a} \gamma + i \Lambda F_{B+\hat a}\gamma) .
\end{align}

\subsection{Énergie}
\label{subenergie}
Soit $(A,\Phi)$, où $\Phi=(\alpha,\gamma)$ et $A=B+a$,  une
configuration à support dans $t > t_0$ telle que $\Phi-\Phi_0$ et
$a$ aient une décroissance en $O(e^{-\epsilon t})$ ainsi que toutes
leurs dérivées ; cette hypothèse
 sert à justifier les intégrations par parties.   On
calcule alors d'après les formules (\ref{Weitz1}) et
(\ref{Weitz2}) en intégrant par rapport à la métrique $\gmod_\tau$
\begin{align*}
  \int |\nabla_{\hat a}\alpha|^2  &= \int 2|\delbar_{\hat a}\alpha|^2 +
  \langle i\Lambda F_{\hat a} \alpha,\alpha\rangle = \int 2|\delbar_{\hat a}\alpha|^2 +
 i\Lambda F_{\hat a} |\alpha|^2 \\
  \int |\nabla_{\hat a}\gamma|^2 &= \int 2|\delbar_{\hat a} ^*\gamma|^2 -
  i\Lambda F_{B+ \hat a} |\gamma|^2 \;;
\end{align*}
en additionnant ces identités, il vient
\begin{align*}
  \int \frac 12 |\nabla _A \Phi |^2 & = \int  \frac 12 (|d_{\hat a}
  \alpha|^2 + |\nabla_{\hat a} \gamma|^2) \\
 &=  
\int  |\delbar_{\hat a} \alpha|^2 +  |\delbar^*_{\hat a} \gamma|^2 +
 \frac i2\Lambda F_{\hat a} |\alpha|^2
 - \frac i2 \Lambda F_{B+\hat a} |\gamma|^2  \\
&=  \int  |\delbar_{\hat
 a}\alpha  +  \delbar^*_{\hat a}\gamma|^2 - 2\langle \delbar_{\hat
 a}\alpha  ,  \delbar^*_{\hat a}\gamma \rangle   +
 \frac i2 ( \Lambda F_{\hat a} |\alpha|^2
 -  \Lambda F_{B+\hat a} |\gamma|^2) \\
&=  \int  |\delbar_{\hat
 a}\alpha  +  \delbar^*_{\hat a}\gamma|^2 - 2\langle \delbar ^2_{\hat
 a}\alpha  ,  \gamma \rangle   +
 \frac i2 ( \Lambda F_{\hat a} |\alpha|^2
 -  \Lambda F_{B+\hat a} |\gamma|^2) 
\end{align*}
or $ \delbar^2_{\hat a}= F_{\hat a} ^{0,2} + N\partial_{\hat a}$ et $2F_{\hat a} = F_a$, donc
\begin{multline*}
\int |\delbar_{\hat
 a}\alpha  +  \delbar^*_{\hat a}\gamma|^2 = 
\int \frac 12 \left [ |d_{\hat a}
  \alpha|^2 + |\nabla_{\hat a} \gamma|^2 \right ] \\
 +\int  \frac i2\Lambda F_B|\gamma|^2 - \frac i4 \Lambda
  F_a(|\alpha ^2|-|\gamma |^2 +\shyp) + \langle
  F_a ^{0,2}\alpha,\gamma\rangle+ 2\langle
  N\partial_a \alpha,\gamma\rangle + \shyp\frac i4  \Lambda F_a.
\end{multline*}
En outre
$$  \langle F_a ^{0,2}\alpha,\gamma \rangle =   \langle F_a ^{0,2},\bar \alpha
  \gamma \rangle = - \left |F_a ^{0,2} - \frac {\bar\alpha\gamma}2
  \right |^2 +|F_a
  ^{0,2} | ^2 +  \left |\frac {\bar\alpha\gamma}2 \right |^2 
$$
et 
$$
i \Lambda
  F_a(|\alpha ^2|-|\gamma |^2 +\shyp) =  \left |i\Lambda F_a
  + \frac {|\alpha|^2-|\gamma|^2 +\shyp}2 \right|^2 -\left |i\Lambda
  F_a\right |^2- \left |\frac
  {|\alpha|^2-|\gamma|^2 +\shyp}2 \right | ^2.
$$
On remarque que
\begin{equation}
\label{idcourb}
  \left | \frac {|\alpha |^2 - |\gamma|^2 +\shyp}4 \right |^2 + 
\left  |\frac {\bar \alpha \gamma}2 \right |^2  =  \left | \frac {|\alpha |^2 +
    |\gamma|^2 +\shyp}4 \right |^2 -  \frac
\shyp 4 {|\gamma|^2}  
\end{equation}
et on en déduit l'identité
\begin{multline}
\label{formulekm}
\int   |\delbar_{\hat
 a}\alpha  +  \delbar^*_{\hat a}\gamma|^2 
+\frac 14  \left |i\Lambda F_a
  + \frac {|\alpha|^2-|\gamma|^2 +\shyp}2 \right|^2
+  \left |F_a ^{0,2} - \frac {\bar\alpha\gamma}2
  \right |^2  - \shyp\frac
 i4  \Lambda F_a
=  \\
\int \frac 12 \left [ |d_{\hat a}
  \alpha|^2 + |\nabla_{\hat a} \gamma|^2 \right ] + \left [ \frac
 i2\Lambda F_B - \frac \shyp 4 \right ]|\gamma|^2   + \frac 1{16} \left (
  |\alpha|^2 +|\gamma|^2 +\shyp \right ) ^2 \\ 
 + \int  2\langle
  N\partial_a \alpha,\gamma\rangle + |F_a ^{0,2}|^2 + \frac 14 \left |i\Lambda
  F_a\right |^2.
\end{multline}
Puisque $a$ et $\Lambda F_a$ décroissent par hypothèse  en
$O(e^{-\epsilon t})$, on en déduit en intégrant par parties que 
  $\int_{t\geq T} \Lambda F_a = \int_{t=T} \frac i4
a\wedge \omega$. Comme par hypothèse $\omega$ s'étend en une forme
fermée sur la variété $M$ tout entière, on en déduit pour tout
$\epsilon>0$ une constante $C$
indépendante de $\tau$ telle que 
\begin{equation}
\label{majorekm}
 \left |\int _{t\geq T} \shyp \frac i4 \Lambda F_a \vol\right | \leq  C+
\epsilon \|da\|_{L^2(\gmod_\tau, t<T)}^2. 
\end{equation}

Suivant \cite{KM} on introduit alors une \emph{énergie} <<à la
Taubes>>~:
\begin{equation}
E_{T}^\tau(A,\Phi)= \int_{t\geq T} 
\left  [ |\nabla_A (\alpha,\gamma)|^2 + |F_a^+ |^2 +|\gamma|^2   
+ \left (
  |\alpha|^2 + |\gamma|^2 +\shyp \right ) ^2\right ]  \vol ^{g_\tau}
\end{equation}
où les normes et les connexions sont prises relativement à la
métrique $g_\tau$.\medskip

Nous allons montrer que l'énergie des solutions des équations de Seiberg--Witten
est uniformément bornée dans le cas où $g$ est suffisamment proche de
$\gmod$ ce qui est automatiquement vrai, à l'action près d'un
difféomorphisme, dans le cas d'une métrique d'Einstein
asymptotique à la métrique de Bergmann.  Dans le cas général d'une
métrique d'Einstein asymptotiquement hyperbolique complexe, il faut légèrement modifier
la construction des métriques $\gmod$ et $\gmod_\tau$
afin que l'hypothèse de la proposition \ref{propenergie} suivante
continue d'être vérifiée (cf. section \ref{secbootstrapping}).

 Cette propriété d'énergie bornée est essentielle pour
faire converger  une suite de solutions pour les métriques $g_\tau$
vers une solution des équations de Seiberg--Witten pour la métrique limite $g$.

\begin{prop}
\label{propenergie}
Supposons que $g$ soit une métrique  telle que  $g=\gmod + O(e^{-\delta
t})$ avec $\delta >2$. Alors, il existe une constante $C$ telle que  
pour $T$ suffisamment grand, 
l'énergie  des solutions $(A_\tau,\Phi_\tau)$ des équations de
Seiberg--Witten (\ref{SWp1},\ref{SWp2}) vérifie   $E ^\tau_{T}
(A_\tau, \Phi_\tau)\leq C$ pour tout $\tau \geq T$.
\end{prop}
\begin{demo}
Soit $(A,\Psi)$ une solution des équations de Seiberg--Witten
perturbées pour une
des métriques $g_\tau$ dans une jauge où $A\otimes B_\tau^{-1}$ et
$\Psi -\Phi_0$ décroissent exponentiellement vite (cf. \cite{KM}
cor. 3.16). En multipliant $\Psi$ par la fonction de 
troncature  $\chi_T$, on peut appliquer l'identité (\ref{formulekm}) à  
$(A,\Phi) = (A,\chi_T\Psi)$ ;
pour $\kappa$ aussi petit que l'on veut, on sait d'après le lemme
\ref{courbchern} que $  \frac i2\Lambda
F_{B_\tau} - \frac \shyp 4  \geq - \kappa$ pour $t \geq T$
suffisamment grand. Par
ailleurs, pour $t\geq T$ assez grand, on peut également supposer que
$|N|\leq\kappa$. Ainsi, le membre de droite de 
(\ref{formulekm}) est un majorant de
\begin{equation}
\label{marredeca}
c \int _ {t\geq T} |\nabla_A (\alpha,\gamma)|^2 + | d a ^+ |^2 -\frac
\kappa c
|\gamma |^2 + (|\alpha|^2 +|\gamma|^2 +\shyp)^2 \vol^{\gmod_\tau},
\end{equation}
où $c>0$ est une constante bien choisie. En utilisant l'identité
(\ref{formulekm}) et l'inégalité
(\ref{majorekm}), on en déduit que (\ref{marredeca}) est majoré par 
\begin{equation}
\label{marredeca2}
 C +  \epsilon \|da\|^2 + \int _ {t\geq T}   |\delbar_{\hat
 a}\alpha  +  \delbar^*_{\hat a}\gamma|^2 
+\frac 14  \left |i\Lambda F_a
  + \frac {|\alpha|^2-|\gamma|^2 +\shyp}2 \right|^2
+  \left |F_a ^{0,2} - \frac {\bar\alpha\gamma}2
  \right |^2 \vol^{\gmod_\tau}.
\end{equation}
Par hypothèse $g-\gmod = O(e^{-\delta t})$ avec $\delta > 2$ et les
connexions spinorielles associées vérifient  
\begin{equation}
\label{connexions}
|\nabla^{\gmod_\tau}_A -
\nabla_A ^{g_\tau}| = O (e^{-\delta t}(1-\chi_\tau(t)).
\end{equation}
 Comme $e^{-\delta
t}\in L^2(g)$, on en
déduit  en utilisant la borne  $C^0$ sur le spineur que pour des constantes $T, C_2>0$ suffisamment grandes,
(\ref{marredeca}) est minoré par 
\begin{equation}
\label{marredeca3}
\frac c2 \int _ {t\geq T} |\nabla_A (\alpha,\gamma)|^2 + | d a ^+ |^2   -
\epsilon |da|^2 -\frac{2\kappa}c
|\gamma |^2 + (|\alpha|^2 +|\gamma|^2 +\shyp)^2 \vol^{g_\tau} - C_2 ,
\end{equation}
où les normes, la connexion et le projecteur $+$ sont pris par rapport à la métrique
$g_\tau$. En appliquant les équations de Seiberg--Witten au terme
$d^+a$ avec l'identité (\ref{idcourb}), il est facile de voir, à
condition d'avoir choisi $\kappa$ suffisamment petit qu'il existe une
constante $c'>0$ telle que (\ref{marredeca3}) soir minoré par
\begin{equation}
\label{marredeca4}
c' \int _ {t\geq T} |\nabla_A (\alpha,\gamma)|^2 + | d a ^+ |^2   -
\frac{c\epsilon}{c'} |da|^2 + |\gamma |^2 + (|\alpha|^2 +|\gamma|^2 +\shyp)^2 \vol^{g_\tau} - C_2 ,
\end{equation}

En utilisant les équations de Seiberg--Witten, la propriété
(\ref{connexions}) et la borne $C^0$ sur le spineur 
on en déduit que l'expression (\ref{marredeca2}) est
majorée par 
$$C + C_2 + 2\epsilon \|da\|^2_{L^2(g_\tau,M)} =  C + C_2 + 4\epsilon
\|da^{+_{g_\tau}}\|_{L^2(g_\tau,M) }^2 ,
$$
quitte à choisir $T$ suffisamment grand ;
c'est donc également un majorant de (\ref{marredeca4}). Nous insistons
sur le fait que $c$ et $c'$ ne dépendent pas du choix de
$\epsilon$. 

En résumé,
nous avons  montré que pour tout $\epsilon' >0$, il existe des
constantes $T,C'>0$,  telles que  
$$ C  + \epsilon
\|da_\tau ^{+_{g_\tau}}\|_{L^2(g_\tau,M) }^2 \geq E_T^\tau(A_\tau, \Phi_\tau).
$$
En découpant
$$\epsilon' \|da_\tau ^{+_{g_\tau}}\|_{L^2(g_\tau,M) }^2 = \epsilon'\|da_\tau
^{+_{g_\tau}}\|_{L^2(g_\tau,t< T) }^2 + \epsilon'\|da_\tau
^{+_{g_\tau}}\|_{L^2(g_\tau,t>T) }^2  
$$
on remarque que le premier terme est borné à l'aide de la borne
uniforme $C^0$ sur le spineur, et que le deuxième terme est contrôlé
par l'énergie en choisissant $\epsilon'$ suffisamment petit. On en déduit la proposition.
 \end{demo}\medskip

\subsection{Première valeur propre du laplacien}
La courbure moyenne est un outil essentiel pour analyser des
différents laplaciens que nous rencontrerons par la suite.  

Soit un produit riemannien tordu $[t_1,t_2]\times Y$ muni d'une métrique,
 $g=dt^2+g_t$, où $g_t$ est une famille de métriques sur $Y$. On définit la courbure  moyenne des sphères $S_t=\{t\}\times Y$ de rayon $t$ par 
$$H= \frac 13 \trace_g \ffond = -  \frac 13 \frac{\dt \vol^{g_t}}{\vol^{g_t}},
$$
où $\ffond$ désigne la deuxième forme fondamentale de $S_t$
et $\vol^{g_t}$ est la forme volume associée à $g_t$. Par
commodité, on définit également une courbure moyenne normalisée
$h=\frac 32 H$.
On a alors le lemme très utile suivant :
\begin{lemme}
\label{courburemoyenne}
Soit $g=dt^2+g_t$ une métrique riemannienne sur $[t_1,t_2]\times Y$ et 
$\delta\in\R$. 
Supposons qu'il existe $h_0\in \R$
 tel que 
$0 > h_0 \geq h$. Alors pour toute fonction $f$, on a 
\begin{equation}
  \int_{[t_1,t_2]\times N}  \!\!\!\! \!\!\!\! |\dt f|^2
 \vol^g 
 \geq h_0^2  \int_{[t_1,t_2]\times N}   \!\!\!\! |f|^2\vol^{g_t} 
+h_0 \Bigg [ \int_{t=t_2 }  \!\!\!\!\!  |f|^2\vol^{g_t}
- \int_{t=t_1}\!\!\!\!\!  |f|^2\vol^g \Bigg ]. 
\end{equation}
\end{lemme}
En particulier  la métrique asymptotiquement hyperbolique $\gmod$
vérifie $\vol ^\gmod =  \sh^2t \sh (2t) dt\wedge\eta\wedge d\eta$ d'où
$h \rightarrow-2$ ; on peut lui appliquer le lemme précédent. En
revanche, pour une métrique $\gmod_\tau$, la courbure moyenne tend vers
$0$ près de l'infini sur le cône asymptotiquement plat. Néanmoins on a
le contrôle suivant :
\begin{lemme}
\label{lemmenormetau}
Pour
toute fonction $f$ à support compact définie sur 
$ M \cap \{t\geq  \tau\}$,  on a
\begin{equation}
 \int _{t>\tau }|\dt f|^2 \vol ^{\gmod_\tau}  -  (\tanh\tau)^{-2}\int_{\tau } |f|^2 \vol ^{\gmod_{\tau,t}}
 \geq  \int_{t>\tau} \left | \frac {\ch\tau}{\xi_2} f \right | ^2
 \vol ^{\gmod_\tau}. 
\end{equation}
les normes étant prises par rapport à la métrique $\gmod_\tau$.
\end{lemme}
\begin{demo}
On fait les changements
de variables $\xi_2= \ch \tau (t-\tau) + \sh \tau $ puis $ r = \ln
\xi_2 $. Alors $\frac {d\xi_2}{\xi_2} = dr$, $d\xi_2 = \ch\tau dt$ et la forme
volume de 
$\gmod_\tau$ s'écrit $ 2 \xi^3_2d\xi_2 \wedge\eta\wedge
d\eta$ ; on en déduit l'identité 
$$\int _{t>\tau}|\dt f|^2\vol^{\gmod_\tau} =  \int 2 
\ch^2 \tau| \partial_{\xi_2}  f|^2  \xi_2^3   d\xi_2 \wedge\eta\wedge
d\eta = \int 2\ch^2\tau | \partial_r  f|^2 e^{2r}    dr \wedge\eta\wedge
d\eta . 
$$
En appliquant le lemme \ref{courburemoyenne}  avec $h_0=-1$, le terme de bord sur
$t_2$ est nul près de l'infini car $f$ est supposée à support compact,
puis en faisant le changement de variable inverse, on obtient
\begin{equation}
\int_{t>\tau} |\dt f|^2 \vol^{\gmod_\tau} \geq \int_{t>\tau}\left |\frac {\ch
\tau}{\xi_2} f\right |^2 \vol^{\gmod} + (\tanh\tau) ^{-2} \int_{t = \tau} | f|^2
\vol^{\gmod _{\tau,t}}
\end{equation}
d'où le résultat.
\end{demo}

 On 
introduire naturellement un nouvel espace de Sobolev : pour toute  fonction
$f$  sur $M$, on pose
\begin{equation}
\label{normetau}
 \| f \|_{\tau} = \|f\|_{L^2(g_\tau, t<\tau)}  +   \left \|\frac
{\ch \tau}
{\xi_2} f\right \|_{L^2(g_\tau, t> \tau)}.
\end{equation}

\begin{lemme}
\label{lemmenormetau2}
Il existe une constante $c>0$ et $T$ suffisamment grand tels que pour
toute fonction $f$ à support compact définie sur 
$Z_T = M \cap \{t\geq  T\}$ nulle en $t=\tau$ et pour toute métrique
$g_\tau$ avec $\tau> T$  on ait
\begin{equation}
\|\dt f\|^2_{L ^2(g_\tau, Z_T)}  \geq c \|  f \|_{\tau, Z_T} ^2 ;
\end{equation}
ce lemme est également vérifié par les métriques $\gmod_\tau$.
\end{lemme}
\begin{demo}
On commence par démontrer le lemme pour les métriques $\gmod_\tau$,
puis on en déduit automatiquement le résultat pour les métriques
$g_\tau$. Il suffit de constater que les intégrants sont perturbés
lorsqu'on passe de la métrique $\gmod_\tau$ à la métrique $g_\tau$ par
un terme  de la forme $(1-\chi_\tau(t))\epsilon(t) (|f|^2+|d
f|^2)$, où $\chi_\tau$ est la fonction de troncature définie au
$(\ref{chitau})$
 et $\epsilon(t)$ une fonction indépendante de
$\tau$ tendant vers $0$. Ce terme provient du recollement entre
$g$ et $\gmod$ pour $ t<\tau$ et
$\epsilon(t)\rightarrow 0$
tient compte  du fait que $g$ est asymptotique à $\gmod$. 

Mais ce résultat est évident pour les métriques $\gmod_\tau$ : on
choisit une constante $0<c\leq 1$, on minore
$$ \int _ {[T,+\infty]} |\dt f |^2\vol  \geq \int _ {[T,\tau]} |\dt f
|^2\vol +  c \int _ {[T,+\infty]} |\dt f |^2\vol,
$$
et on ajuste $c$ de sorte qu'après avoir appliqué les lemmes
\ref{courburemoyenne} et \ref{lemmenormetau} il ne reste plus qu'une
intégrale de bord positive en $t=\tau$ dans le membre de droite de l'inégalité.
\end{demo}

On en déduit classiquement une inégalité de Poincaré sur $M$.
\begin{cor}[Inégalité de Poincaré]
\label{poincare}
Il existe une constante $c>0$ telle que pour toute fonction
 $f$  sur $M$ telle que $\|f\|_\tau<\infty$, on ait
$$ \|df\|_{L^2(g_\tau)} \geq c \|f\|_\tau.
$$
\end{cor}
\begin{demo}
Supposons l'assertion fausse. On en déduit une suite de fonctions
$f_\tau$ telles que $ \|df\|_{L^2(g_\tau)} \rightarrow 0$ et
$\|f\|_\tau =1$. Alors, quitte à extraire, cette suite admet une
limite faible sur tout 
compact  $f\in L^2_1(g)$. A la limite $df=0$ donc $f$ est
constante. Comme la métrique est de volume infini et $f$ est $L^2$,
nécessairement $f=0$. Par compacité de l'inclusion $L^2\subset
L^2_1$, la suite $f_\tau$ converge en fait  vers $0$ au sens $L^2$-fort sur tout compact
de $M$. 

Choisissons $T$ suffisamment grand afin d'appliquer le
lemme~\ref{courburemoyenne2} ;
on majore 
\begin{align*}
 1 & = \| f_\tau \|_\tau ^2  =  2\| (1-\chi_T)  f _\tau \|^2 +
2\|\chi_T f_\tau \|_ \tau ^2\\
&\leq   2\| (1-\chi_T) f_\tau \|^2 + 2c^{-1} \int |d( \chi_T  f_\tau)|^2 \vol \\
&\leq   2\| (1-\chi_T) f_\tau \|^2 + 4c^{-1} \int |(\dt\chi_T) f_\tau
|^2\vol + 2c^{-1} \int |\chi_T d f_\tau|^2  \vol, 
\end{align*}
mais le membre de droite tend vers $0$ lorsque $\tau$ tend vers l'infini
d'où une contradiction. Le corollaire est par conséquent vérifié.
\end{demo}\medskip

\subsection{Laplacien des $1$-formes}
\label{calclapl}
Calculons le laplacien des $1$-formes pour une
métrique $g= dt + \xi_1^2\eta^2 +\xi_2^2\gamma$, telle que $\xi_1
=2\xi_2\xi_2'$. Une $1$-forme $a$ se décompose en 
$$a= fdt+q\xi_1\eta + b,
$$
où $f$ et $q$ sont des fonctions et $b$ une $1$-forme tangente à la
distribution de contact, c'est à dire telle que $i_\dt b =i_\XX b
=0$. Un calcul direct de $\Delta a= (dd^* + d^*d)a$, en utilisant la
formule $\dt b = \nabla_\dt b- \ffond (b)$,  montre que
\begin{align*}
\Delta a &= \left (-\partial_t^2 f +2h \dt f +2h'f + \Delta f -
x^2\cdot f  +\frac {\xi'_1}{\xi_1} x \cdot q + [\dt,
d^*]b  \right )dt \\
 &+ \left (-\partial_t^2q  +2h \dt q +2h'q + \Delta q -
x^2\cdot q  -\frac {\xi'_1}{\xi_1} x \cdot f  + [x,
d^*]b  - 2\frac {\xi_2'}{\xi_2} *db \right )  \xi_1\eta \\
 &-\nabla_\dt^2 b  +2h \nabla_\dt b - \frac {(\xi_2')^2 +
2\xi_2\xi_2''}{\xi_2^2} b + 2\Delta b  + *x*x\cdot b - 2\frac {\xi_2'}{\xi_2}(df - *dg )  ,
\end{align*}
où $x=\xi_1^{-1}R$, avec $R$ le champ de Reeb,  et les opérateurs $d, * ,d^*,\Delta$ apparaissant
dans le membre de droite sont définis le long de la distribution de
contact. Plus précisément, $d b$ représente la projection
orthogonale de la
différentielle de $b$ sur la distribution de contact ; $*$ est
l'opérateur de Hodge défini sur la distribution
de contact  par $b\wedge *b = \xi_2^2d\eta$ ; on en déduit une
divergence  $d^*$ et un laplacien $\Delta$ formés à l'aide de ces deux opérateurs.

En intégrant par parties pour $t$ fixé et en utilisant le fait que la
forme volume est invariante suivant $x$, on obtient
\begin{align}
\label{lpl1tranche}
\int_t \langle \Delta a,a \rangle \vol^{g_t}  &= \int _ t \Big
(\langle -\partial_t^2 f +2h \dt f, f \rangle +
|x\cdot f |^2 + \langle \frac {\xi'_1}{\xi_1} x \cdot q , f\rangle  +
|d f|^2 +  
2h'|f|^2  \\
 &+  \langle-\partial_t^2q  +2h \dt q, q \rangle   +
|x\cdot q| ^2  - \langle \frac {\xi'_1}{\xi_1} x \cdot f ,g\rangle +
|d q|^2 + 
2h'|q|^2   \nonumber \\
 &+ \langle -\nabla_\dt^2 b  +2h \nabla_\dt b - \frac {(\xi_2')^2 +
2\xi_2\xi_2''}{\xi_2^2} b, b \rangle  + 
|x\cdot b |^2 + |d b|^2 +  |d^* b|^2  \Big )\vol^{g_t}  .\nonumber 
\end{align}
En intégrant maintenant par parties suivant $t$, on obtient
\begin{align*}
\int_{[t_1, t_2]} \langle \Delta a,a \rangle \vol^{g} & + \int_{\partial[t_1,t_2]}
\langle \nabla_\dt a,a \rangle \vol^{g_t} \\
 = \int _{[t_1,t_2]} &\Big (
|\dt f|^2 +
|x\cdot f |^2 + \langle \frac {\xi'_1}{\xi_1} x \cdot q , f\rangle  +
|d f|^2 +  
2h'|f|^2  \\
 &+  |\partial_t q|^2   +
|x\cdot q| ^2  - \langle \frac {\xi'_1}{\xi_1} x \cdot f ,g\rangle +
|d q|^2 + 
2h'|q|^2   \\
 &+ |\nabla_\dt b|^2  - \frac {(\xi_2')^2 +
2\xi_2\xi_2''}{\xi_2^2}| b|^2  + 
|x\cdot b |^2 + |d b|^2 +  |d^* b|^2  \Big )\vol^{g}  .
\end{align*}
Par la formule de Bochner $\Delta a = \nabla ^*\nabla a + \Ric(a)$, on
reconnait que le terme de gauche de l'identité ci-dessus est égal à
$\int _ {[t_1,t_2]} (|\nabla a|^2 + \Ric(a,a) )\vol^g$.

Par construction, la courbure moyenne des métriques $\gmod_\tau$ est
croissante, soit $h'\geq 
0$. En minorant,  
$$  \langle \frac {\xi'_1}{\xi_1} x \cdot q , f\rangle
+ |x\cdot f|^2 -  \langle \frac {\xi'_1}{\xi_1} x \cdot f , q\rangle
+ |x\cdot q|^2 \geq -  |\frac {\xi'_1}{2\xi_1}  f |^2 -  
|\frac {\xi'_1}{2\xi_1} q |^2,$$
 on en déduit l'inégalité
\begin{multline}
\label{contrlpl1}
\int_{[t_1,t_2]} \left (|\nabla a|^2 +\Ric (a,a)\right ) \vol^{g} \geq \\
 \int _{[t_1,t_2]} \left (
|\dt f|^2   -  |\frac {\xi'_1}{2\xi_1} f |^2 
 +  |\partial_t q|^2   -  |\frac {\xi'_1}{2\xi_1}  q |^2
+ |\nabla_\dt b|^2  - \frac {(\xi_2')^2 +
2\xi_2\xi_2''}{\xi_2^2}| b|^2   \right )\vol^{g}  .
\end{multline}
Ce calcul nous permet de démontrer le lemme suivant.
\begin{lemme}
\label{courburemoyenne2}
Il existe des constantes $c,T>0$ telles que pour toute métrique
$g_\tau$ et toute
$1$-forme $a\in L^2_1(g_\tau)$ définie sur le bout  $[T,+\infty[\times
Y$ on ait 
$$\int_{t>T} \left (|\nabla a|^2 + \Ric^{g_\tau}(a,a) \right ) \vol \geq c
\int_{[T,\tau]} |a|^2\vol,
$$
les normes étant prises par rapport à la métrique $g_\tau$. Ce lemme
est également valable pour les métriques $\gmod_\tau$.
\end{lemme}
\begin{demo}
On commence par démontrer le lemme pour les métriques $\gmod_\tau$,
puis on en déduit automatiquement le résultat pour les métriques
$g_\tau$. Il suffit de constater que les intégrants sont perturbés
lorsqu'on passe de la métrique $\gmod_\tau$ à la métrique $g_\tau$ par
un terme  de la forme $(1-\chi_\tau(t))\epsilon(t) (|a|^2+|\nabla
a|^2)$, où $\chi_\tau$ est la fonction de troncature définie au
(\ref{chitau})
et  $\epsilon(t)$ une fonction indépendante de
$\tau$ tendant vers $0$.

Examinons les termes d'ordre $0$ dans le membre de droite de
l'inégalité (\ref{contrlpl1}). Sur l'intervalle $[T,\tau]$, on a par
construction 
$$ \frac {\xi_1'}{2\xi_1} \leq  1 ,\quad  \frac {(\xi_2')^2 + 2 \xi_2
\xi_2''}{\xi_2^2}
\leq 3 ; $$
sur la partie conique  $[\tau,+\infty]$, on a 
$$ \frac {\xi_1'}{2\xi_1} =  \frac {\ch \tau}{2\xi_2} ,\quad  \frac {(\xi_2')^2 + 2 \xi_2
\xi_2''}{\xi_2^2} =  \left (\frac {\ch \tau}{\xi_2}\right )^2.$$
On utilise le lemme \ref{courburemoyenne} avec $h_0$ proche de $-2$
sur l'intervalle $[T,\tau]$ et le lemme \ref{lemmenormetau} sur
l'intervalle $[\tau,+\infty]$. On en déduit, en négligeant le terme de
bord en $T$ qui est positif,  une constante $c>0$ telle que
$$\int_{[T,\tau]} \left (|\nabla a|^2 + \Ric^{g_\tau}(a,a) \right )
\vol \geq  c \int_{[T,\tau]}  | a|^2 \vol +  h_0 \int_{\tau}  | a|^2 \vol^{g_t},$$
 et
$$\int_{[\tau, +\infty]} \left (|\nabla a|^2 + \Ric^{g_\tau}(a,a) \right )
\vol \geq   (\tanh \tau)^{-2} \int_{\tau}  | a|^2 \vol^{g_t}.$$
On peut ajuster  une constante $0 <c'\leq 1 $ telle que $c'h_0 +
(\tanh\tau)^{-2}>0$ pour tout $\tau$ suffisamment grand. Alors
\begin{align*}
  \int _ {[T,+\infty]}\left (|\nabla a|^2 + \Ric^{g_\tau}(a,a)
\right ) \vol   \geq c' &\int_{[T,\tau]} \left (|\nabla a|^2 + \Ric^{g_\tau}(a,a)
\right )\vol \\ + 
& \int_{[\tau,+\infty]}  \left (|\nabla a|^2 + \Ric^{g_\tau}(a,a)
\right )\vol ,
\end{align*}
et le terme de bord en $\tau$, apparaissant dans le membre de droite  après minoration est
 positif ;
  on obtient ainsi 
l'inégalité souhaitée.
\end{demo}

\subsection{Convergence des solutions sur tout compact}
Nous avons maintenant tous les outils nécessaires pour faire converger
la suite de solution  des équations de Seiberg--Witten
$(A_\tau,\Phi_\tau)$ sur tout compact de $M$. 
Nous supposerons dans la fin de cette section que  l'énergie des
solutions   $(A_\tau, \Phi_\tau)$ admet une borne uniforme. 
Nous voulons faire converger cette
suite quitte à faire des changements de jauge et à extraire une
sous-suite vers une solutions des équations pour la métrique  $g$
\begin{eqnarray}
  \label{SW1}
  \Dirac ^ {g} _{A} \Phi & =  & 0  \\
  F ^{+_{g}}_{A\otimes B^{-1}} & = & q(\Phi) - q(\Phi_0) ,\nonumber
\end{eqnarray}  
où $B$ est la connexion de Chern de la métrique $\gmod$.
 
Dans un premier temps nous faisons converger  les solutions hors d'un
compact suffisamment grand. On commence par fixer une jauge de Hodge
en faisant agir le groupe de jauge $\jauge= \Map(M,S^1)$.
\begin{lemme}
\label{lemmehodge}
Quitte a faire des changements de jauge, on peut supposer que 
 $d^*a_\tau = 0$ et $a_\tau \in L^2(g_\tau)$ sur $M$.
\end{lemme}
\begin{demo}
On cherche une fonction  $u$ telle que $d^*a_\tau + d^*du = 0$.
Il est possible de résoudre le problème en minimisant la fonctionnelle
$$ \int_M (\frac 12 |du| +\langle a_\tau , du\rangle )\vol^{g_\tau}
$$
et on obtient une solution telle que
$$\int |du|^2 \vol = -2\int \langle a
,du\rangle \leq 2\|du\|_{L^2(g_\tau)}\|a\|_{L^2(g_\tau)}, $$
d'où $du \in {L^2(g_\tau)}$. 
En fait, on en sait un peut plus sur $u$ car d'après le corollaire \ref{poincare}
$$\|du\|_{L^2(g_\tau)}\geq c \|u\|_\tau
$$
pour toute fonction à support compact et l'inégalité reste vraie par
passage à la limite.
Par conséquent, $u$ admet une borne $L^2_1$ sur chaque compact et
par régularité elliptique locale, on en déduit  que $u$ est une
fonction lisse. Finalement, la transformation
de jauge  $e^{u/2}$ résout le problème.
\end{demo}

Notons qu'on perd a priori de la régularité sur $\Phi_\tau$ en se plaçant en
jauge de Hodge, mais l'énergie bornée offre un contrôle suffisant pour
faire converger les solutions. La borne $C^0$ a priori sur le spineur
et l'énergie fournissent une borne uniforme sur $\|d ^ +a_\tau
\|_{L^2(g_\tau)}$. En intégrant par parties et en utilisant le fait
que nous sommes dans une jauge de Hodge, on en déduit une borne
uniforme sur $\int_M|\nabla a_\tau|^2 + \Ric(a_\tau,a_\tau)
\vol$. Puis en utilisant le lemme \ref{courburemoyenne2}, on en déduit
une borne uniforme sur $\int_{[T, \tau]}|a|^2\vol$ pour $T$
suffisamment grand. On en déduit quitte à extraire une sous-suite que
$a_\tau$ converge sur tout compact vers une limite faible $a\in
L^2_1(g)$. 

La borne $C^0$ sur le spineur  $\Phi_\tau$ et la borne
 sur
$\int_{t>T} |\nabla \Phi_\tau|^2\vol$ fournie par l'énergie nous
indiquent que $\Phi_\tau$ est borné en norme $L^2_1$ sur un compact
fixé. Par extraction diagonale, on en déduit que $\Phi_\tau$ converge
sur tout compact de $M$ vers une limite faible $\Phi$. On obtient ainsi une
 limite $(A,\Phi)$, où $A=B+a$, solution des équations de (\ref{SW1})
 et d'énergie finie,  définie sur le bout  $t>T$. Par ellipticité du
 problème $(A,\Phi)$ est en fait lisse. L'énergie nous permet
 maintenant de récupérer une contrôle global sur $(A,\Phi)$.

\begin{lemme}
\label{lemmeregsup}
Soit $(A,\Phi)$ une solution d'énergie finie des équations de Seiberg--Witten
(\ref{SW1}) pour la métrique $g$ telle que $A\otimes B^{-1}\in
L^2_1(g)$ et $\Phi$ borné en norme $C^0$. Alors quitte à faire un changement de jauge par une
constante $u\in S^1$ on a $\Phi-\Phi_0\in L^2_1$.
\end{lemme}
On peut donc faire converger les solutions $(A_\tau,\Phi_\tau)$
sur le bout de $M$. On se place en jauge de Hodge à l'aide du lemme
\ref{lemmehodge} et  on peut extraire une limite sur tout compact
$(A,\Phi)$ d'énergie finie. Puis du  lemme \ref{lemmeregsup} on déduit :
\begin{cor}
\label{convbout}
Supposons l'énergie des solutions
$(A_\tau,\Phi_\tau)$ uniformément bornée. Alors pour $T$ suffisamment
grand, quitte à faire des
changements de jauge et à extraire une sous-suite,  on en déduit que
$(A_\tau,\Phi_\tau)$ converge sur tout compact de $M\cap\{t\geq T\}$
vers une solutions $(A,\Phi)$  des équations de 
Seiberg--Witten 
(\ref{SW1}), relatives à la métrique asymptotiquement hyperbolique
complexe $g$, telle que $A\otimes B^{-1}$ et $\Phi-\Phi_0$ sont dans $L^2_1(g)$.
\end{cor}
\remarque au départ nous avons défini des métriques $g_\tau$ pour une
famille d'indices $\tau\in \R^+$. Dans le corollaire ci-dessus et dans
toute la suite de cet article, une \emph{suite
extraite} de $(v_\tau)_{\tau\in \R^+}$ désignera une sous-suite
$(v_\tau)_{\tau\in I}$, où $I\subset \R^+$ est un ensemble
infini et non borné.\medskip 

\begin{demode}{du lemme \ref{lemmeregsup} }
L'énergie finie nous permet de conclure instantanément que $\gamma\in
 L^2_1(g)$ en décomposant $\Phi= (\alpha,\gamma)$. Quant à $\alpha$,
 on sait seulement que $\int_{t>T}|\nabla_{\hat a}\alpha | ^2\vol
 <+\infty$. En identifiant  par transport
 parallèle les spineurs sur le bout de $M$ à des
 familles de spineurs sur $Y$ dépendant de $t$, on peut écrire
$$ \Phi (t_2,y) -   \Phi (t_1,y) = \int _{t_1}^{t_2} (\nabla^A_{\dt}\Phi) dt.
$$
On en déduit que 

$$ |\Phi (t_2,y) -   \Phi (t_1,y)|^2 \leq \frac {e^ {-4t}}4 \int _{t_1}^{t_2} |\nabla^A_{\dt}\Phi|^2 e^{4t}dt.
$$
Comme $|\nabla^A\Phi|^2$ est intégrable, on en déduit que $\Phi(t,y)$
est de Cauchy pour presque tout $y$. On obtient une limite $\Phi(y)$
telle que 
\begin{equation}
\label{pourquoi}
 \int _{\{t_1\}\times Y} |\Phi (t,y) -   \Phi (y)|^2 \vol ^g \leq c \int
_{[t_1,+\infty]\times Y} |\nabla^A_{\dt}\Phi|^2 \vol^g.
\end{equation}

Évidemment par construction $\nabla_ \dt^A\Phi(y)=0$, et via le lemme
\ref{courburemoyenne} on obtient 
$$ \int_{[T,t_1]} |\nabla^A_{\dt}\Phi|\vol^g\geq h^2_0 \int_{T,t_1} |\Phi
- \Phi(y)|^2 \vol^g +  h_0 \int_{\partial[T,t_1]} |\Phi
- \Phi(y)|^2 \vol^{g_t} ;
$$
D'après (\ref{pourquoi}) l'intégrale de bord en $t_1$ tend vers $0$
lorsque $t_1$ tend vers l'infini. On en déduit que $\Phi - \Phi(y) \in
L^2(g)$. 

Fixons un point $y\in Y$. 
Puisque la métrique est asymptotiquement hyperbolique complexe, 
on a une suite d'ouverts $V_y(t)$, qui s'identifient à la boule unité
 centrée en $0$ de $\CHH$  via une
isométrie, sur lesquels la métrique $g$ tend vers la métrique
hyperbolique complexe lorsque $t$ tend vers l'infini cf. (\cite{B} section I.1.B). En transportant
également la solution $(A,\Phi)$ des équations de Seiberg--Witten sur
$V_y(t)$, on en déduit une  suite de
solution $(A_t,\Phi_t)$  sur la boule unité de $\CHH$ centrée en
$0$. Puisque $(A,\Phi)$ était d'énergie finie, l'énergie de
$(A_t,\Phi_t)$ tend vers $0$ sur la boule. Comme par hypothèse $A\otimes B^{-1} \in
L^2_1(g)$, on en déduit que $A_t\otimes A_0^{-1}$, où $A_0$ est la
connexion de Chern du fibré anti-canonique de $\CHH$, tend vers $0$
en norme $L^2_1$,
relativement à la métrique hyperbolique complexe, lorsque $t$ tend
vers l'infini. Puisque l'énergie de $(A_t,\Phi_t)$ est bornée et que
$\Phi_t$ est uniformément borné en norme $C^0$, on obtient une borne
uniforme sur la norme $L^2_1$ de $\Phi_t$ sur la boule unité de
$\CHH$. On peut donc en extraire une limite faible $\Psi$. A la
limite, $(A_0,\Psi)$ est une solution d'énergie nulle des équations de
Seiberg--Witten pour la métrique hyperbolique complexe. On en déduit
que $\Psi$ est parallèle et finalement que $\Psi=(u\sqrt{-\shyp},0)$ où
$u\in S^1$ est constante.
Par compacité de l'inclusion $L^2\subset L^2_1$, $\Phi_t$ converge en
fait vers $\Psi$ au sens $L^2$-fort sur la boule. Par
conséquent en ramenant $\Phi(y)$ sur la boule unité, on a $\Psi=
\Phi(y)=(u\sqrt{-\shyp},0)$. En faisant agir la transformations de jauge constante
$u^{-1}$ sur $(A,\Phi)$, on en déduit le lemme.
\end{demode}

En ce qui concerne le domaine compact $M\cap\{t\leq T+1\}$, la borne uniforme $C^0$
sur le spineur suffit pour faire converger la suite de solution
$(A_\tau,\Phi_\tau)$. Nous nous référons par exemple au lemme $4$ de
\cite{KM2} afin d'énoncer le résultat suivant :
\begin{lemme}
\label{convbord}
Soit une variété à bord $\overline X$ munie d'une structure $\spinc$
et d'une métrique riemannienne $g$. Soit $(A_\tau,\Phi_\tau)$ une
suite de solutions des équations de Seiberg--Witten telle qu'il existe
$C>0$ avec $|\Phi_\tau|\leq C$ 
pour tout $\tau$.  
Quitte à faire des changements de jauge et  à
extraire une sous-suite, la suite $(A_\tau,\Phi_\tau)$ converge au
sens $C^\infty$ sur
$\overline X$ vers
une solutions $(A,\Phi)$ des équations de Seiberg-Witten (\ref{SW1}).
\end{lemme}
\remarque c'est la même démonstration que dans le cas compact. Il faut
faire attention  pour faire converger la connexion : sur les variétés
à bord, on dispose d'une théorie de Hodge en fixant une condition de
Neumann convenable pour les formes harmoniques. Cette condition se
traduit sur les $1$-formes par l'annulation suivant la normale au
bord. 

On peut également supposer que la métrique dépend de $\tau$ dans le
lemme \ref{convbord}  et converge au sens $C^\infty$
vers une métrique $g$ sur $\overline X$.

Notre problème est maintenant de <<recoller>> les solutions obtenues
sur les domaines $t\leq T +1 $ et  $t\geq T$ de $M$.
Nous avons besoin du lemme préliminaire suivant :
\begin{lemme}
\label{apriori2}
Soit  $(A,\Phi)$ une solutions   d'énergie finie des équations de Seiberg--Witten
(\ref{SW1}) pour la métrique $g$ , définie près de l'infini. Alors il
existe  $\epsilon  >0$,
et  $T$ suffisamment grand  tels que sur le domaine $M\cap \{t\geq
T\}$ on ait $|\alpha| >\epsilon$, en décomposant $\Phi=(\alpha,\gamma)$.
\end{lemme}
\begin{demo}
Supposons le lemme faux. Alors, il existe une suite de points $z_j\in
M$ tendant vers le bord à l'infini tels que
$|\alpha|(z_j)\rightarrow 0$.  En faisant  la même construction que dans
la démonstration du lemme \ref{lemmeregsup}, on obtient une suite
de  $(A_j,\Phi_j)$ en ramenant $z_j$ en $0$ par une isométrie  de
$\CHH$ dont l'énergie tend vers $0$ sur la boule unité. Quitte à faire des
transformations de jauge, cette suite converge  vers une solutions
d'énergie nulle  des
équations de Seiberg--Witten sur $\CHH$. Par
régularité elliptique, la convergence est en réalité $C^\infty$ et
 donc à la limite $\alpha (0)=0$ ;  ceci contredit l'énergie nulle.
\end{demo}\medskip

En appliquant le lemme \ref{convbout}
 on fait converger la suite 
 $(A_\tau,\Phi_\tau)$  vers $(A,\Phi)$ sur le bout  $t \geq T $ pour
 $T$ suffisamment 
 grand  en extrayant et en faisant des transformations de
 jauge sur $M$. Puis en appliquant le lemme \ref{convbout}, on trouve une
 suite de transformations de jauge $u_\tau$ telles que, quitte à
 extraire, $u_\tau\cdot (A_\tau,\Phi_\tau)$ converge sur $M \setminus
 \{t > T+1\}$. En décomposant $\Phi_\tau=(\alpha _\tau,\gamma_\tau)$ et en utilisant le
 fait que les convergences sont $C^\infty$ sur tout compact, on en
 déduit d'après le lemme \ref{apriori2} un $\epsilon >0$ tel que pour
 $\tau$ suffisamment grand, $|\alpha_\tau|>\epsilon $ pour tout $t\geq
 T$. On en déduit que $u_\tau$ converge également sur l'anneau
 $[T,T+1]$. Alors pour $N$ et $\tau$
 suffisamment grand, on peut écrire $u_\tau = e^{iv_\tau}u_N$  pour un
 unique $v_\tau$ tel que $|v_\tau| < \pi/2 $.  On définit les
 transformations de jauge $f_\tau$ par : $f_\tau= u_\tau u_N^{-1} $ pour $t\leq
 T$, $f_\tau = \exp(i(1-\chi_{T+1})v_\tau )$ pour $t\geq T$
 (cf. (\ref{chitau}) pour la définition de $\chi$) . Alors $f_\tau
 \cdot (A_\tau,\Phi_\tau)$ converge sur $M$ tout entier et on
peut énoncer le résultat  suivant :
\begin{prop}
\label{propexistence}
Soit $g$ une métrique asymptotiquement hyperbolique complexe sur $M$, munie
d'une structure $\spinc$ adaptée. Soit $(A_\tau,\Phi_\tau)$ une suite
de solutions d'énergie bornée des équations de Seiberg--Witten
perturbées relatives aux 
métriques asymptotiquement plates $g_\tau$ approximant $g$.
Alors quitte à faire des changements de jauge et à extraire une
sous-suite, $(A_\tau,\Phi_\tau)$ converge sur tout compact de $M$ vers une
solution $(A,\Phi)$ des équations de Seiberg--Witten (\ref{SW1}) pour la
métrique $g$ telle que $\Phi-\Phi_0 \in L^2_1$,
$a=A\otimes B^{-1}\in L^2_1$ et $d^*a = 0 $ près de l'infini.
\end{prop}
\medskip
	
\section{Métriques d'Einstein asymptotiquement hyperboliques complexes}
\label{secbootstrapping}
\bigskip
On suppose maintenant dans toute cette section que $g$ est une
métrique d'Einstein asymptotiquement hyperbolique complexe. Nous
expliquons maintenant comment faire en sorte que la perturbation des
équations de Seiberg--Witten (\ref{SWp2}) s'évanouisse et comment généraliser
les résultats obtenus pour ces nouvelles constructions.

\subsection{Modification des constructions}
D'après la proposition 3.1 de \cite{BH} (nous nous référons largement
 à sa démonstration dans ce qui suit), il existe une métrique
 de Kähler--Einstein asymptotiquement hyperbolique complexe $(\bar
 g,\bar J)$ d'infini conforme $\gamma$ définie sur le bout de
$M$ par une série formelle. Nous
 supposerons pour simplifier 
dans la suite de cette section comme si la métrique était
 exactement Kähler--Einstein plutôt que définie par une série formelle
 bien que nos résultats restent vrais  en toute généralité. 

Précisons quelques points sur la construction de $(\bar g,\bar J)$.
 La structure complexe $\bar
 J$ est obtenue en modifiant $J$ le long de la distribution de
 contact :
la différence entre $J$ et $\bar J$ est paramétrisée par un tenseur
 $\theta \in H^{0,1}\otimes H_{1,0}$ avec $H^{0,1}$ les formes de type
 $(0,1)$ et $H_{1,0}$ les vecteurs de type $(1,0)$ de la distribution
 de contact tel que
$$ T_{0,1}^{\bar J} =\{ X + \theta (X), X\in  T_{0,1}^{J} \}.
$$
 De plus on sait que par construction $\theta = O(e^{-2t})$.
Les tenseurs $\theta_\tau= (1 - \chi_{\tau-1})\theta$ définissent
une famille de structures presque complexes qui réalisent un recollement entre 
$\bar J$ pour $t\leq \tau-2$ et $J$ pour $t\geq \tau -1$. On en déduit
une famille $\tilde J_\tau$ de structures presque complexes en
remplaçant $J$ par $J_\tau$ pour $t\geq \tau-1$.

La métrique presque kählérienne $\gmod = dt^2 + \sh^2(2t)\eta^2+
\sh^2(t)   \gamma$ découle, comme dans le cas hyperbolique complexe
(\ref{potentiel}), 
d'un potentiel $U= -\frac 12 \ln (1-\tanh^2 t )$ avec $\omega^{\gmod}=
i\partial\delbar U$. La métrique $\bar g$ est définie en modifiant ce
potentiel par $\omega^{\bar g}= i\partial\delbar (U+F)$ où $F=O(e^{-2t})$.
De façon analogue au cas de la structure presque complexe, la suite de
potentiels  $F_\tau= (1 - \chi_{\tau-1})F$, permet de définir une suite de
métriques presque kählériennes $\tilde g_\tau$ relativement à $\tilde
J_\tau$ recollant la métrique $\bar g$ pour $t\leq \tau-2$ et
$\gmod_\tau$ pour $t\geq \tau -1$.

Toutes nos constructions gardent un sens en remplaçant la suite $(\gmod_\tau, J_\tau)$
de la section \ref{subhc}
par $(\tilde g_\tau, \tilde J_\tau)$.  On définit alors la suite
d'approximations $g_\tau$ de $g$ par (\ref{defgtau}). 
Comme $\bar g$ et $\gmod$  diffèrent par un $O(e^{-2t})$,  on en déduit que l'ensemble des résultats de
la section \ref{secbulle}, le corollaire \ref{poincare} et le lemme
\ref{courburemoyenne2}  restent vrais pour les nouvelles métriques
$g_\tau$. Il s'en suit que les solutions de Seiberg--Witten perturbées
associées aux métriques $\gmod_\tau$ admettent une uniforme borne
$C^ 0$ suivant le lemme \ref{controleC0}.

Si $g$ est une métrique  d'Einstein   asymptotiquement hyperbolique
complexe, elle satisfait une équation elliptique modulo
l'action des difféomorphismes. Un examen précis de cette équation mené
par Biquard et Herzlich donne
des renseignements  supplémentaires sur le
comportement à l'infini de $g$~:
\begin{resultat}[Biquard--Herzlich]
Si $g$ est une métrique d'Einstein asymptotiquement hyperbolique,
alors quitte à faire agir un difféomorphisme, on a près de l'infini
\begin{equation}
\label{exprbh}
 g= \bar g + k + O(e ^{-(4+\epsilon)t}), \mbox{ où $\epsilon>0$}
\end{equation}
 et  $k= e^{-2t}k_0$  avec $k_0$ une $2$-forme bilinéaire symétrique sur la
distribution de contact de $Y$, anticommutant à la structure
complexe. Autrement dit $i_{\eta^\sharp} k_0=0$ et $k_0(J\cdot,J\cdot) = -k_0$.
\end{resultat}

Comme la métrique $\bar g$ est de Kähler--Einstein, l'équation
$F_B^{+_{\bar g}} = q(\Phi_0)$ est vérifiée.  On en déduit que
$F_B^{+_g} = q(\Phi_0)+ \varpi$, avec $\varpi= O(e^{-4t})$.
Les équations de Seiberg--Witten perturbées sont maintenant
définie par 
\begin{eqnarray*}
\Dirac_{A_\tau}^{g_\tau} \Phi_\tau &=&0\\
F^+_{A_\tau} - q(\Phi_\tau) &=& F^+_{ B_\tau }  - q(\Phi_0) + \varpi_\tau
\end{eqnarray*}
où $B_\tau$
est la connexion de Chern induite par la 
métrique $(\tilde g_\tau, \tilde J_\tau)$ sur le fibré anti-canonique
et $\varpi_\tau= (1-\chi_\tau) \varpi$.
La perturbation $F^+_{ B_\tau }  - q(\Phi_0) + \varpi_\tau $ est
maintenant nulle  pour $t\leq \tau-2$. À la limite on obtient les 
équations de Seiberg--Witten non perturbées (\ref{SWnp}) pour la
métrique $g$.

Avec les mêmes démonstrations que celles des propositions
\ref{propenergie}, \ref{propexistence}, on obtient en utilisant le
résultat de Biquard--Herzlich
\begin{prop}
Soit $g$ une métrique d'Einstein asymptotiquement
hyperbolique complexe et
soit $(A_\tau,\Phi_\tau)$ une suite de solutions des équations de
Seiberg--Witten perturbées pour les approximations $g_\tau$.
Alors pour $T$ suffisamment grand, 
l'énergie  $E ^\tau_{T} (A_\tau, \Phi_\tau)$ est uniformément bornée.
\end{prop}
Comme la borne uniforme sur l'énergie est systématique, on a donc :
\begin{prop}
\label{propexistence2}
Soit $g$ une métrique d'Einstein asymptotiquement hyperbolique complexe sur $M$, munie
d'une structure $\spinc$ adaptée. Soit $(A_\tau,\Phi_\tau)$ une suite
 des équations de Seiberg--Witten
perturbées relatives aux 
approximations  asymptotiquement plates $g_\tau$  de  $g$.
Quitte à faire des changements de jauge et à extraire une
sous-suite, $(A_\tau,\Phi_\tau)$ converge sur tout compact de $M$ vers une
solution $(A,\Phi)$ des équations de Seiberg--Witten non perturbées
(\ref{SWnp}) pour la 
métrique $g$ telle que $\Phi-\Phi_0 \in L^2_1$,
$a=A\otimes B^{-1}\in L^2_1$ et $d^*a = 0 $ près de l'infini, avec
$(A_0,\Phi_0)$  la solution standard des équations pour une
métrique formelle Kähler--Einstein $\bar g$, asymptotiquement
hyperbolique complexe, définie sur le bout de $M$.
\end{prop}
\medskip

Si $g$ est d'Einstein, nous allons montrer, quitte à faire
un changement de jauge, que la
solution $(A,\Phi)$  possède une forte  décroissance
à l'infini en nous servant de~(\ref{exprbh}).\medskip

\subsection{Décroissance des solutions}
Par commodité, on considère maintenant des solutions $(A,\Phi)$ des
équations de Seiberg-Witten, où $A$ est par convention une
connexion définie sur la racine carrée du fibré déterminant
$\C^{1/2}$. Nous retrouvons ainsi une solution des équations
précédentes par un simple changement de variables
$(A,\Phi)\rightarrow (A \otimes A, \sqrt 2\Phi)$ et l
'action du groupe de jauge $\jauge = \Map (M,S^1)$  est donnée par 
$$u\cdot(A,\Phi) = (A +\frac {du}u,u ^{-1}\Phi).
$$

En linéarisant les équations de Seiberg--Witten et l'action du
groupe de jauge $\mathcal G$ en la solution standard  $(B,\Phi_0)$ des équations
pour la métrique de Kähler Einstein $\bar g$, on
obtient les opérateurs
\begin{align*}
\delta_1 (u) &= (du,-u\Phi_0) \\
\delta_2 (a,\phi) &= (d^+a-\{\Phi_0\otimes \phi ^* +\phi\otimes
\Phi_0^*\}_0, \Dirac_B\phi + a\cdot \Phi ).
\end{align*}
L'adjoint formel de $\delta_1$ est donné par
$$\delta_1^* (a,\phi) = d^*a - i\Im(\Phi_0^*\phi). $$
On définit l'opérateur $\mathcal D=\delta_1 ^*+\delta_2$
 et son  laplacien $P = \mathcal D^*\mathcal D$.
En notant l'espace de Sobolev à poids $L^2_{k,\delta} = e^{-\delta t}L^2_k$, on a la proposition~:
\begin{prop}
\label{isomsw}
Il existe des constantes $\epsilon, c, T>0$   telles
        que pour tout couple $(a,\phi)$ à support compact définit
        sur $[T,+\infty [$ nul en  $t=T$ et
        tout $\delta \in [0,2+\epsilon]$ on 
        ait
\begin{equation}
\int |\mathcal D (a,\phi)|^2 \poids\vol  \geq c \int \left
(|\nabla a |^2+ |\nabla_{A_0}\phi|^2+|a|^2+|\phi|^2\right)\poids\vol,
\end{equation}
où $\mathcal D$ est l'opérateur des équations linéarisées modulo
l'action du groupe de jauge en
en $(B,\Phi_0)$. 
\end{prop}\medskip

Il n'est
pas difficile d'en déduire le corollaire suivant en notant 
$$\|\sigma\|^2_{L^2_ {k,\delta}(g) } =  \sum_{j=0}^k \|(\nabla^k
\sigma) e^{\delta t}\|^2_{L^2(g) }.
$$

\begin{cor}
\label{isomswcor}
Il existe des constantes $\epsilon,T, c_k >0$   telles
        que pour  toute configuration $(a,\phi)$ à support compact
        nulle en $t = T$   et    tout $\delta \in [0,2+\epsilon]$ on 
        ait
\begin{equation}
\| \mathcal D (a,\phi) \|_{L^2_{k-1,\delta }(g)}  \geq c_k \| (a,\phi)
\|_{L^2_{k,\delta }(g)}  . 
\end{equation}
\end{cor}\medskip

Un procédé de bootstrapping, dans le cas où $g$ est une métrique d'Einstein, va nous permettre en s'appuyant sur ce
corollaire d'obtenir une borne $L^2_k$ sur une
solution $(A,\Phi)$ des équations de Seiberg--Witten non perturbées
(\ref{SWnp}) en partant d'une borne $L^2_1$.

Si la métrique vérifie $g= \bar g +
O(e^{-(4+\epsilon)t} )$, on en déduit automatiquement que  
\begin{equation}
\label{vardirac}
\Dirac_A^{g } \Phi = \Dirac_A^{ \bar g} \Phi +
 O(e^{-(4+\epsilon)t})\cdot \Phi,
\end{equation}
où le  $O(e^{-(4+\epsilon)t})$  est la $1-$forme donnée par
 $\nabla_A^{g}- \nabla_A^{\bar g}$ et agissant sur $\Phi$ par produit
 de Clifford. En utilisant le fait que $\Phi$ est borné en norme
 $C^0$ on en déduit $ O(e^{-(4+\epsilon)t})\cdot \Phi$ est aussi un $ O(e^{-(4+\epsilon)t})$.

 Mais d'après (\ref{exprbh}), la métrique
$g$ comporte également une perturbation donnée par un tenseur
 $k=O(e^{-4t})$ défini sur la distribution de contact et
anticommutant avec $J$. 
La variation de l'opérateur de
Dirac en fonction de la métrique est calculée  précisément
dans \cite{BG} au théorème $21$ : en identifiant les fibrés de
spineurs par un isomorphisme canonique à celui réalisé par la métrique
$\bar g$, la variation infinitésimale de l'opérateur de Dirac pour une
variation $ k$ de 
métrique $\bar g$ est donnée par
$$ \dot \Dirac^{\bar g} _A  \Phi = - \frac 12 Cl(\nabla_{ K}^{A,
\bar g} \Phi) +\frac 14 (\delta ^{\bar g} k + d\trace_{\bar g} k)\cdot \Phi,
$$
où $K$ est l'endomorphisme de $TM$ défini par $\bar g^{-1}k$ et $Cl$
est la contraction par le produit de Clifford.
On applique cette identité à $B$ et $\Phi_0$ et puisque $\Phi_0$ est
parallèle le premier terme est nul. Comme $k$ anticommute à $J$ sa
trace par rapport à $\gmod$ est nulle. 
Par ailleurs on a  $\delta^\gmod k = O(e^{-5t})$
(cf. formule $(7.16)$ de \cite{BH} en faisant attention à la
convention différente sur la courbure scalaire). Comme $\gmod - \bar
g= O(e^{-2t})$, on en  déduit que
pour un accroissement  $k=O(e^{-4t})$, on a 
$\delta^{\bar g} k + d\trace_{\bar g} k = O(e^{-5t})$ d'où
$$ \dot \Dirac^{\bar g} _{B}  \Phi_0 = O(e^{-5t}).
$$
Les termes quadratiques et d'ordre supérieur  dans la variation de l'opérateur de Dirac par
rapport à la métrique  $(\Dirac^g_A - \Dirac^{\bar g}_A)\Phi$ ont une
décroissance encore plus forte. On en déduit finalement que 
la formule~(\ref{vardirac}) reste valable pour toute
métrique d'Einstein asymptotiquement hyperbolique complexe. 

Notons que le tenseur $k$ correspond à une  déformation
infinitésimale de la
métrique $\bar g$ laissant la forme symplectique invariante. Il est
facile d'en déduire que l'équation de Seiberg--Witten portant sur la courbure est
perturbée  par un terme à décroissance en
$O(e^{-(4+\epsilon)t})$ lorsqu'on passe de $\bar g$ à $g$.

Rappelons que $(A,\Phi)$ est dans une jauge de Hodge près de l'infini
et les équations de Seiberg-Witten
s'écrivent  par conséquent sous la forme
\begin{equation}
\label{bootstr}
\mathcal D(a,\phi) + Q(a,\phi) + O\left
(e^{-(4+\epsilon)t} \right )= (i\Im\langle\Phi_0, \phi \rangle,0,0)
\end{equation}
où $\phi = \Phi-\Phi_0$, $a= A\otimes B^{-1}$, $Q$ est un terme quadratique et $\mathcal D$ est l'opérateur des
équations de Seiberg--Witten pour la métrique $\bar g$ linéarisées en
$(B,\Phi_0)$. Le terme
quadratique $\mathcal Q(a,\phi)$ admet une borne $L^2_1$. En effet 
$$\mathcal Q(a,\phi)= (0, a\cdot \phi , q(\phi,\phi) ) ;
$$
maintenant $\nabla q(\phi,\phi) = 2q(\nabla \phi, \phi) \in L^2$ puisque $\nabla
 \phi \in L^2$ et $\phi$ borné en norme $ C^0$ et de même $(\nabla
 a)\cdot\phi \in L^2$ car $\nabla a \in L^2$. Il nous reste à voir que
 $a\cdot\nabla \phi \in L^2$. Par l'inclusion de Sobolev $L^2_1\subset
 L^4$, on a $a\in L^4$. Par ailleurs, on a d'après les équations de
 Seiberg--Witten $\Dirac_{B}\phi = - a\cdot\phi$. Comme  $\phi \in
 C^0$, on en déduit  que  $  \Dirac_{B} \phi \in L^4$. De plus 
$\phi\in L^4$ et par ellipticité de l'opérateur de Dirac, on en déduit une borne $L^4_1$ sur $\phi$. Finalement
$\nabla \phi\in L^4$ et on a une  borne $L^2$ sur $\nabla \mathcal
 Q(a,\phi)$.
 
\remarque  puisque les opérateurs  $\Dirac_{B}$
est défini à partir d'une métrique asymptotiquement hyperbolique complexe,
les constantes de régularité elliptique locale peuvent être choisies
indépendamment du centre des boules ce qui justifie le raisonnement de
bootstrapping ci-dessus. 

D'après (\ref{bootstr}) on a donc une borne $L^2_1$ sur $\mathcal
D(a,\phi)$. On en déduit une borne $L^2_1$ sur  $\mathcal D(\chi_T (a,\phi))$, où
$\chi_T$ est la fonction de troncature définie à la section
\ref{subhc}, puis via le corollaire \ref{isomswcor} on en déduit que
$(a,\phi)$ est dans $L^2_2$. En réitérant ce raisonnement, on montre
que $(a,\phi)\in L^2_k$ pour tout $k$.

\begin{cor}
La solution $(A,\Phi)$ des équations de Seiberg--Witten non perturbées
(\ref{SWnp}) pour la métrique d'Einstein asymptotiquement hyperbolique
complexe $g$, en jauge de Hodge près de l'infini, obtenue dans la
proposition \ref{propexistence2}
vérifie $A\otimes B ^{-1}\in
L^2_k$ et  $\Phi - \Phi_0\in
L^2_k$ pour tout $k$.
\end{cor}

On  a maintenant obtenu une régularité suffisante sur la solution $(A,\Phi)$
pour  fixer une jauge de Coulomb près de l'infini.
\begin{lemme}
\label{lemmecoulomb}
Soit $(A,\Phi)$ une configuration telle que $A-B\in L^2_k$ et
$\Phi-\Phi_0\in L^2_k$ avec $k\geq 2$. Alors on peut faire un changement  de jauge
sur le bout de $M \cap \{t>T\}$, avec $T$ suffisamment grand, de la
forme $e^{u}$, avec $e^u-1\in L^2_{k+1}$ tel qu'on ait
 $\delta^*_1(A-B,\Phi-\Phi_0)=0$, où $\delta_1$ est
l'action linéarisée du groupe de jauge en $(B,\Phi_0)$.
\end{lemme}
\begin{demo}
 On introduit  les espaces $W_k$ qui sont les
complétés pour la norme $L^2_k$ des fonctions (formes, sections)  à
support compact
définies pour $t\geq T$ et nulles en $t=T$. On considère l'action du groupe de
jauge  entre les espaces 
\begin{eqnarray*}
\{  (A,\Phi) / a,\phi\in W_k  \mbox{ et } \delta_1^*(a,\phi)=
0\}\times \{ u \in W_{k+1} \} & \longrightarrow&   \{
(A,\Phi) / a,\phi\in W_k  \} \\
((A,\Phi), u ) &\longmapsto & (A + du , e^{-u}\Phi),
\end{eqnarray*}
où l'on a noté $(a,\phi) = (A\otimes B^{-1},\Phi-\Phi_0)$. Cette
flèche est définie à condition de choisir  $k\geq 3$, afin d'avoir
l'inclusion continue $L^2_k \subset C^0$.
Nous allons voir que c'est un difféomorphisme local près de $((B,\Phi_0),0)$.
Il suffit de montrer que la différentielle en ce point est un
isomorphisme, c'est à dire de montrer qu'étant donné $(\dot
a,\dot\phi)\in W_k$ il existe un unique $ u \in W_{k+1}$ qui
résout infinitésimalement le problème de la jauge de Coulomb 
$ \delta_1^*\delta_1( u ) = - \delta ^*_1 (\dot a, \dot\phi)$. Il suffit
donc de montrer que $\delta_1 ^*\delta_1 : W_{k+1} \rightarrow
W_{k-1}$ est un isomorphisme. On calcule aisément
$$\delta ^*_1 \delta_1 (u)= d^*du +|\Phi_0|^ 2 u $$
d'où
$$ \int \langle \delta^* _1 \delta_1 (u), u \rangle \vol \geq
|\Phi_0|^2 \int |u|^2 \vol ;
$$
et on en déduit une constante $c>0$
 telle que  
\begin{equation}
\label{coulombboot}
\|u\|_{L^2}\leq c\|\delta ^*_1 \delta_1 (u)\|.
\end{equation}
Par ailleurs, compte tenu du fait que la métrique est asymptotiquement
hyperbolique, les constantes $c_k$ des estimations elliptiques locales sur
les boules de rayon $1$
$$\|u\|_{L^2_{k+2}}\leq c_k\left ( \|u\|_{L^2} +\|\delta _1^*
\delta_1 (u)\|_{L^2_k} \right )$$ 
peuvent être choisie indépendamment du centre de la boule. En partant
de (\ref{coulombboot}), on en déduit des constantes $c'_k$ telles que 
$$\|u\|_{L^2_{k+2}}\leq c'_k\|\delta _1^*
\delta_1 (u)\|_{L^2_k},$$
où les normes sont prises sur $M$ tout entier ; l'opérateur
$\delta_1^* \delta_1$ est donc un isomorphisme.

Soit $(A,\Phi)$ une configuration telle que dans l'énoncé du
corollaire. Notons $(a,\phi)= (A\otimes B^{-1},\Phi-\Phi_0)$
et posons
$A_T =  B + \chi_T a$ et $\Phi_T = \Phi_0 +
\chi_T \phi$ où $\chi_T$ est la
fonction de troncature définie à la fin de la section
\ref{subhc}. Lorsque $T$ tend vers l'infini, $(A_T,\Phi_T)$ tend vers 
$(B,\Phi_0)$ en norme $L^2_k$. On en déduit la jauge de Coulomb
souhaitée pour $T$ suffisamment grand.  
\end{demo}
\medskip

\remarque soit $\sigma$ une section de fibré (forme, spineur) telle
que  
$\sigma \in L^2_{k,\delta} $. Alors
$$\sigma(t)= -\int_t^\infty e^{-(4+2\delta )t}(\nabla_\dt\sigma)
e^{(4+2\delta)t} dt ;
 $$
par inégalité de Cauchy-Schwarz, il vient
\begin{equation}
 |\sigma(t)|^2 \leq \frac {e^{-2(2+\delta)t}}{4+2\delta}
\int_t^\infty |\nabla_\dt \sigma|^2 e ^{(4+2\delta )t}dt,\label{rqdecr}
\end{equation}
En utilisant l'inclusion continue $L^2_2\subset C^0$ sur $Y^3$, on en
déduit en intégrant sur $Y$ que $\sigma,\nabla \sigma,\dots,\nabla^{k-3}\sigma$
sont des $O(e^{-(2+\delta) t})$

  En particulier,
on sait automatiquement que pour la solution des équations de
Seiberg-Witten obtenue dans la proposition~\ref{propexistence},
placée dans une jauge de Coulomb près de l'infini via le lemme
\ref{lemmecoulomb}, 
$A-B$ et $\Phi-\Phi_0$ et leurs dérivées ont une décroissance en $O(e ^ {-2t})$.

Du corollaire  \ref{isomswcor}, nous déduisons maintenant le résultat souhaité pour les
métriques d'Einstein asymptotiquement hyperbolique.
\begin{cor}
\label{cor15}
Soit $(A,\Phi)= (B+a,\Phi_0+\phi)$ une solution des équations de
Seiberg-Witten pour une métrique d'Einstein $g$ asymptotiquement hyperbolique
complexe avec $a,\phi\in L_1 ^2$.  Alors
il existe un $\epsilon >0$ suffisamment petit tel que  dans une jauge
de Coulomb près de l'infini, $a$, $\phi$ et leurs dérivés
d'ordre arbitrairement grand ont une décroissance en $O\left
(e^{-(4+\epsilon)t}\right)$. 
\end{cor}
\begin{demo}
On peut quitte à multiplier par une fonction de troncature
supposer que  $a,\phi$ est nulle en $t=T$. 
Plaçons nous dans une jauge telle que  $\delta_1^*(a,\phi)=
0 $ près de l'infini suivant lemme \ref{lemmecoulomb}.
En reprenant la discussion menant à (\ref{bootstr}), on voit que les
équations de Seiberg-Witten  s'écrivent près de l'infini en jauge
de Coulomb 
\begin{equation}
\label{bootstr2}
\mathcal D(a,\phi) + Q(a,\phi) + O\left
(e^{-(4+\epsilon)t}\right )= 0.
\end{equation}
A partir de cette identité, un argument de bootstrapping basé sur
le corollaire \ref{isomswcor} et la remarque (\ref {rqdecr}) va nous donner la
décroissance voulue :
on sait déjà que  $a,\phi \in L^2_k$ et leurs dérivées jusqu'à
l'ordre $k-3$  ont une décroissance  en $e^{-2t}$ ;
on en déduit en utilisant l'équation~(\ref{bootstr2}) que $\mathcal
D(a,\phi) = O(e^{-4t})$ (où le $O$ porte sur les dérivées
jusqu'à l'ordre $k-3$  ). On a  donc  $\mathcal
D(a,\phi) \in L^2_{k-3, 7/4 }$ et on en déduit en utilisant  le
corollaire~\ref{isomswcor} que $(a,\phi)$ est lui aussi dans
$L^2_{k-3,7/4}$ ; c'est donc un $O(e^{-7/4 t})$ à l'ordre $k-6$ par
la remarque (\ref{rqdecr}). 
On réitère le raisonnement, précédent :  $Q(a,\phi)$ est un
$O(e^{-7/2t})$ \` l'ordre $k-6$ et en choisissant une constante
$0<\epsilon'\leq\min(\epsilon, 8 )$ on en déduit d'après (\ref{bootstr})
que $D(a,\phi) = O(e^{-(4+\epsilon')t})$. On applique à nouveau le
corollaire~\ref{isomswcor} avec $4<\delta<4+\epsilon'$ et on en déduit
que que $(a,\phi) = O(e^{-(4+\epsilon')t})$ à l'ordre $k-6$, d'où
le résultat en partant de $k$ suffisamment grand.
\end{demo}\medskip

On déduit  immédiatement de la
proposition~\ref{propexistence2} et du corollaire~\ref{cor15} le
théorème suivant :
\begin{theo}
\label{theoconv}
Soit $g$ une métrique d'Einstein asymptotiquement hyperbolique et
$g_\tau$ sa suite d'approximations asymptotiquement plates sur une
variété orientable $M^4$ munie d'une structure $\spinc$ adaptée.
Soit $(A_\tau,\Phi_\tau)$ une suite de solutions des équations de
Seiberg-Witten perturbées~(\ref{SWp1},\ref{SWp2}) associées aux 
métriques $g_\tau$. 
Quitte à faire des changements de jauge et à
extraire une sous suite, on peut supposer que $(A_\tau, \Phi_\tau)$
converge sur tout compact au sens $C^\infty$ vers une solution
$(A,\Phi)$ des
équations de Seiberg--Witten non perturbées, 
\begin{eqnarray*}
\Dirac_{A}^{g} \Phi&=&0\\
F^+_{A} &= &q(\Phi) ,
\end{eqnarray*}
telle que 
$A\otimes B^{-1} $, $\Phi-\Phi_0$ et leurs dérivées sont des
$O\left(e^{-(4+\epsilon)t}\right )$ avec $\epsilon>0$.
\end{theo}
\medskip

Le théorème  \ref{theoexistence} est une reformulation  de
ce dernier résultat sans la technique.

\subsection{Laplacien de Seiberg--Witten}
La fin de cette section est dédiée à la démonstration de la
proposition~\ref{isomsw}. 
Nous commençons pas examiner précisément l'opérateur $P$ dans un lemme calculatoire.
\begin{lemme}
\label{calculatoire}
Pour une métrique hermitienne et des champs $(a,\phi)$ nuls en $t=T$,
le laplacien $P$ des équations de Seiberg--Witten  linéarisées en
$(B,\Phi_0)$ relativement à  la métrique de Kähler--Einstein
$\overline g$ vérifie 
\begin{multline}
\int_{t>T} \langle P (a,\phi),(a,\phi)\rangle \poids\vol^{\bar g} = \\
\int_{t>T} \Big\{ \langle\Delta a,a\rangle + \langle\Delta \alpha,\alpha\rangle  + \langle\Delta \gamma,\gamma\rangle  -\frac s2
(|a|^2 +|\alpha |^2 + |\gamma|^2)  \Big\}\poids\vol^{\bar g}.
\end{multline}
\end{lemme}
\begin{demo}
Pour ce calcul, nous aurons besoin des résultats suivants ou
on utilise la convention 
$\pherm AB =\frac 12\trace_\C A^*B$ pour le produit hermitien sur
$\mathrm{End}(W)$~:
\begin{lemme}
\label{lemmeformule}
Étant données $\Phi$, $T$, $\eta$ des sections de $W^+$, $\mathrm{End}_0(W^+)$,
 $\Lambda^{2,+}X$, on a~:
\begin{gather}
|\rho(\eta )|^2 =2|\eta|^2 \label{formule1} \\
 |\sanstr{\spinendi\Phi}|^2 = \frac14 |\Phi|^4 \label{formule2} \\
\pherm T {\sanstr{\spinendii\Phi\phi}} =\frac 12 \pherm {T\Phi}\phi +
\frac 12 \pherm {T\phi}\Phi \label{formule3} \\
|\Phi|^2|\phi|^2 =|\sanstr{\spinendii\Phi\phi}|^2+|\Im \pherm\Phi\phi
|^2 . \label{formule4} 
\end{gather}
\end{lemme}
\bigskip

Nous noterons dans ce qui suit  $f =i_\dt a$. On a alors  la formule de Leibniz 
$$\mathcal D [ ( a,\phi)\poids ]= \poids\mathcal D (a,\phi)
+2\delta\poids(-f,\rho(dt\wedge a)^+, dt\cdot\phi),$$
donc
\begin{equation}
\label{exprPD}
\int \langle P(a,\phi),(a,\phi)\rangle\poids\vol 
=\int  |\mathcal D(a,\phi)|^2 + 2\delta \langle\mathcal D (a,\phi),(-f,\rho(dt\wedge
a)^+,dt\cdot\phi)  \rangle \poids\vol.
\end{equation}
Calculons le premier terme~:
\begin{multline*}
|\mathcal D(a,\phi)|^2  = 
|d ^*a|^2+|\Im\langle\Phi,\phi\rangle|^2 -2\langle
d ^*a,i\Im\langle\Phi, \phi\rangle \rangle \\
+2|d ^+a|^2+|\sanstr{\spinendii\Phi\phi}|^2- 2\langle
d ^+a, \sanstr{\spinendii\Phi\phi}\rangle \\
+|\Dirac_A\phi|^2+|a|^2|\Phi|^2+2\langle\Dirac_A\phi,a\cdot\Phi
\rangle.
\end{multline*}
D'après  (\ref{formule3}) et  (\ref{formule4})
il vient
\begin{multline*}
|\mathcal D(a,\phi)|^2  =  |d ^* a|^2+ 2|d  ^+a|^2
+|\Dirac_A\phi|^2 +2\pherm{\Dirac_A\phi}{a\cdot\Phi} \\
-2\langle(d ^*a+d ^+a)\cdot \Phi,\phi\rangle
+|\Phi|^2(|a|^2+|\phi|^2).
\end{multline*}

Maintenant, intégrons par parties le terme $\pherm
{\Dirac_A\phi}{a\cdot\Phi}$ ;
$$\int \pherm{\Dirac_A\phi}{a\cdot\Phi}\poids\vol =\int 2\delta
\pherm\phi{\dt\cdot a\cdot\Phi}\poids
+\pherm\phi{\Dirac_Aa\cdot\Phi}\poids \vol,$$
or 
\begin{equation}
\Dirac_A (a\cdot\Phi) =
(d  a+d ^*a)\cdot\Phi-a\cdot\Dirac_A \Phi -
2\,\trace \; a \!\otimes \!\!\nabla^A\Phi.
\end{equation}
En effet~:
choisissons une base locale $(e_i)$ orthonormale parallèle en un
point ; alors 
\begin{equation*}
\begin{split}
\Dirac_A ( a\cdot\Phi) & = \sum_{i,j} e_i\cdot\nabla^A_{e_i}
(a_j e_j\cdot\Phi) = \sum_{i,j} d  a_j(e_i) e_i\cdot e_j\cdot\Phi +
a_j e_i \cdot \nabla_{e_i}(e_j) \cdot\Phi+
a_j e_i\cdot e_j\nabla^A_{e_i}\Phi \\
&=(d  a+d ^*a)\cdot\Phi-a\cdot\Dirac_A \Phi - 2\sum_i
a_i \nabla_{e_i}^A\Phi, 
\end{split} 
\end{equation*}
d'où le résultat.

Notons que $\pherm {a\otimes\phi}{\nabla^A\Phi}=\pherm \phi{\trace
\; a \!\otimes \!\!\nabla^A\Phi}$ ; alors
\begin{multline}
\label{exprD}
\int |\mathcal{D}(a,\phi)|^2\poids\vol=\int
|d ^*a|^2+2|d ^+a |^2 +|\Dirac_A\phi|^2 
+|\Phi|^2(|a|^2+|\phi|^2) \\
+4\delta\pherm\phi{dt\cdot a\cdot\Phi}- 
2\pherm\phi{a\cdot\Dirac_A\Phi} - 4 \pherm
{a\otimes\phi}{\nabla^A\Phi} \poids\vol .
\end{multline}
Dans le cas où
$(A,\Phi)=(B,\Phi_0)$, on a 
$|\Phi|^2=-\frac s2$, $\Dirac_A\Phi=0$ et $\nabla^A\Phi=0$~; les deux derniers termes de
l'intégrale sont donc nuls.

Maintenant,
\begin{multline*}
 2\delta \int \langle\mathcal D (a,\phi),(-f,\rho(\dt\wedge
a)^+,dt\cdot\phi)  \rangle \poids\vol = 
2\delta\int \pherm
{d ^*a}{-f}+\pherm f{i\Im\pherm\Phi\phi} +\pherm
{\Dirac_A\phi}{dt\cdot\phi} \\ 
+\pherm{a\cdot\Phi}{dt\cdot\phi}+ 
\pherm{\rho^+(\dt\wedge a)}{\rho^+(d  a)-\sanstr{\spinendii\Phi\phi}
}  \poids\vol 
\end{multline*}
ce qui d'après le lemme~(\ref{lemmeformule}) est égal à 
\begin{equation}
  2\delta\int \pherm{d ^*a}{-f}+2\pherm
{(dt\wedge\gamma)^+}{d ^+a }+2\pherm{a\cdot\Phi}{dt\cdot
\phi} +\pherm{\Dirac_A\phi}{dt\cdot\phi} \poids\vol.
\end{equation}

On applique la formule de Lichnerowicz~:
\begin{equation}
\label{exprlichn}
\int|\Dirac_A\phi|^2\poids\vol = \int
\langle\nabla_A^*\nabla_A\phi,\phi\rangle + \frac s4|\phi|^2 
+\pherm
{F_A^+\cdot\phi}\phi + 2\delta\pherm{dt\cdot\Dirac_A\phi}\phi
\poids\vol .
\end{equation}
En décomposant $\phi=(\alpha,\gamma)$, $\Phi=(\sqrt {-\frac s2},
0)$ il vient
$$\langle F_A^+\cdot\phi,\phi\rangle = \langle q(\Phi)\phi,\phi\rangle
= \frac s4 (|\gamma |^2 - |\alpha |^2).
$$
Finalement, en intégrant par parties on  voit apparaître
\begin{equation}
\label{laplacien2}
\int \langle P (a,\phi),(a,\phi)\rangle \poids\vol = 
\int \Big\{ \langle\Delta a,a\rangle + \langle\Delta \alpha,\alpha\rangle  + \langle\Delta \gamma,\gamma\rangle  -\frac s2
(|a|^2 +|\alpha |^2 + |\gamma|^2)  \Big\}\poids\vol,
\end{equation}
où $\Delta $ désigne le laplacien  de la métrique métrique de
Kähler--Einstein $\bar g$ sur les formes.
\end{demo}

Il s'agit maintenant d'analyser chacun des $3$ laplaciens intervenant
 dans le lemme précédent. Comme la métrique $\bar g$ n'a pas de
 formule explicite, nous allons plutôt faire les calculs dans un
 premier temps  par rapport à la
 métrique asymptotiquement hyperbolique complexe $\gmod = dt^2 +
\sh^2(2t) \eta + \sh^2(t) \gamma$, et de la structure presque complexe
 $J$.

Nous avons déjà traité la question des $1$-formes dans un
 cadre plus général : on part  de (\ref{lpl1tranche}) appliqué  à la 
métrique $\gmod$ et une $1$-forme à support compact
$a=fdt+q \sh(2t)\eta + b$, où $i_\dt b = i_R  b = 0$, 
nulle en $t=T$. Dans ce cas, $\frac {\xi'_1}{2\xi_1} = 1 + O(e^{-t})$,
    $\frac {(\xi'_2)
^2+ 2\xi_2\xi_2''}{\xi_2^2} = 3 + O(e^{-t})$  et on obtient
\begin{align}
\label{lp1}
\int_{t} \langle\Delta^{\gmod} a, a\rangle  \poids\vol^{\gmod_t} \geq 
 \int _{t} \Big \{
 & \Big \langle (-\partial_t ^2+2h\dt -1 + O(e^{-t})  ) f,f \Big
 \rangle \\
+ & \Big \langle (-\partial_t
 ^2+2h\dt -1 + O(e^{-t})) q,q  \Big \rangle \\
+ & \Big \langle (-\nabla_{\partial _ t} ^2+2h\nabla_{\partial _ t} -
 3 +  O(e^{-t})) b,b \Big \rangle 
  \Big \} \poids\vol^{\gmod_t} .
\end{align}

Soit $\Delta^h$  le laplacien défini à l'aide de la différentielle
restreinte aux directions de la structure de contact et de l'opérateur
de Hodge sur la structure de contact donné par $b\wedge *b = \sh^2(t)d\eta$.
Soit le champ de vecteurs $x= \sh^{-1}(2t)R$.
Un calcul direct montre que le laplacien sur les fonctions est donné
sur le bout de $M$ par 
$$ \Delta^{\gmod}\alpha = (-\partial_t^2 + 2h - x^2 + \Delta^h)\alpha.
$$
On en déduit en intégrant par parties sur $Y$ que
\begin{equation}
\label{lp2}
\int_{t} \langle\Delta^{\gmod} \alpha, \alpha \rangle  \poids\vol^{\gmod} \geq 
 \int _{t}  \Big \langle (-\partial_t ^2+2h\dt) \alpha, \alpha \Big
 \rangle \poids\vol^{\gmod} .
\end{equation}

Il faut maintenant s'occuper du laplacien sur les $(0,2)$ formes. Pour
la métrique $\bar g$, le laplacien sur les $(0,2)$ formes est donné
via les identités  kählériennes par $\Delta
\gamma = 2\delbar\delbar^*\gamma$ ; nous allons donc plutôt nous intéresser à
ce laplacien  de Dolbeaut par rapport à la métrique $\gmod$.
On définit le champ de vecteur et la $1$-forme différentielle de types $(1,0)$
$$\XX=\frac  {\dt - i x}2  \quad\mbox { et }\quad \Omega = dt +i
\sh(2t)\eta.  $$
On vérifie facilement que 
\begin{equation}
\label{domega}
d\Omega=\delbar \Omega =  \frac{2i}{\tanh (2t)}dt\wedge \sh
(2t) \eta + \frac {2i}{\tanh(t)} \sh^ 2 (t)d\eta= 2i\omega  ^{\gmod} + O(e^{-t}) ,
\end{equation}
où $\omega^{\gmod}$ est la forme de Kähler de $(\gmod,J)$. D'autre part $\partial\Omega = 0$.

Toute $(0,2)$-forme peut s'écrire de la manière suivante
$$\gamma =  \overline\Omega \wedge \mu,$$
où $\mu$ est une famille de  $(0,1)$-formes sur la distribution de contact
 Maintenant, l'identité kählérienne
$\Delta\gamma = 2\delbar\delbar^*\gamma = -2i
\delbar\Lambda\partial\gamma$ permet de calculer facilement le
laplacien : tout d'abord $\partial \gamma = \partial \overline\Omega
\wedge \mu - \overline\Omega \wedge \partial \mu $~; en notant $\delcp$ les dérivations dans les directions de la distribution de contact, on a
 $\partial \mu = \Omega \wedge \XX\cdot \mu
+\delcp\mu$. Par la formule~(\ref{domega}), il vient
$$\partial \gamma = -2idt\wedge\sh(2t)\eta \wedge ( \XX + \tanh ^{-1}
(2t) )\cdot\mu - (\Lambda\delcp \mu) \overline\Omega\wedge \sh^2(t)d\eta.
$$
 En contractant avec la forme de Kähler $\omega = dt \wedge
 \sh(2t)\eta +\sh^ 2(t) d \eta$, on obtient
$$\Lambda \partial \gamma = - \Lambda\delcp \mu
\overline\Omega -2i(\XX + \tanh^ {-1}(2t))\cdot \mu.
$$
Alors
$$ \Delta ^{\gmod} \gamma= -2i \delbar \Lambda \partial \gamma =
\overline\Omega \wedge \left [\Delta ^h \mu - 
4 \overline \XX \cdot \XX\cdot \mu - 4\overline \XX\cdot
\left (\tanh^ {-1} (2t)\mu\right )\right ].              
$$
Or 
$$4\overline \XX\cdot\XX = \dt ^ 2 + x^ 2
+\frac {2i}{\tanh(2t)} x.$$
Comme $\dt \cdot \mu = \nabla_\dt \mu +  \tanh ^{-1}(t)\mu$, le
laplacien s'écrit finalement
\begin{equation}
 \Delta  ^{\gmod} \gamma = \overline\Omega \wedge \Big [
 -  \nabla_\dt ^ 2 + 2 h \nabla_\dt  
-x ^ 2 
-4i x + \Delta^h -3 + 
O\big (e^{-t}\big)\Big ]\cdot \mu,
\end{equation}
où $\Delta^h\mu = -2i \delbar^h  (\delbar^h)^* \mu$.
Puis en intégrant par parties suivant $Y$, et en majorant
$$ |x\cdot\mu|^2 - \langle 4ix\cdot\mu,\mu\rangle \geq  - 4|\mu|^2
$$

on en déduit que 
\begin{equation}
\label{lp3}
\int_{t} \langle\Delta^{\gmod} \gamma, \gamma \rangle  \poids\vol^{\gmod_t} \geq 
 \int _{t}  \Big \langle (-\nabla_{\dt} ^2+2h\nabla_\dt -7 + O(e^{-t})
 )\gamma, \gamma  \Big
 \rangle \poids\vol^{\gmod} .
\end{equation}
\bigskip

\subsection{Conclusion}
Dans chacun des cas étudiés, le laplacien agit sur une forme $\sigma$ 
de sorte que,
$$\int_t \langle (\Delta^{\gmod}-\frac s2) \sigma,\sigma\rangle
\poids\vol^{\gmod} \geq \int_t  \langle
(-\nabla_\dt ^2 \sigma +2h\nabla_\dt \sigma + \kappa) \sigma  ,
\sigma\rangle\poids \vol^{\gmod}.
$$
où d'après (\ref{lp1},\ref{lp2},\ref{lp3}), $\kappa$ est  une constante
 telle que $\kappa \geq  -s/2-7 +O(e^{-t}) = 5   +O(e^{-t}) $.
 En
intégrant par parties suivant $t$ et en utilisant l'hypothèse de
 nullité de $\sigma$ en $t=T$,  il vient 
$$
\int_{t\geq T} \langle
-\nabla_\dt ^2 \sigma ,
\sigma\rangle +2h\langle\nabla_\dt  \sigma ,
\sigma\rangle\poids \vol^{\gmod}   = \int 2\delta  \langle
\nabla_\dt\sigma,\sigma \rangle + |\nabla_\dt  \sigma|^2 \poids \vol^{\gmod}.
$$
Par ailleurs, pour  $c>0$ 
$$-2\langle \nabla_\dt\sigma,\sigma\rangle \leq c ^{-1} |\nabla_\dt
 \sigma| ^2+ c|\sigma| ^2.$$
En appliquant les deux inégalités précédentes, il vient
\begin{align*}
\int \langle (\Delta^{\gmod}-\frac s2 )\sigma, \sigma\rangle \poids\vol^{\gmod} & \geq
\int 2(\delta - c) \langle \nabla_\dt\sigma,\sigma\rangle + (\kappa
- c^2)|\sigma |^2 \poids \vol^{\gmod} \\
& \geq \int \left [2(\delta - c)(h - \delta) + (\kappa
- c^2)\right ]|\sigma |^2 \poids \vol^{\gmod}.
\end{align*}
En choisissant $c=\delta-h$,
on obtient l'inégalité
\begin{equation}
\label{isom00}
\int_{t\geq T} \langle (\Delta^{\gmod}-\frac s2 ) \sigma,\sigma\rangle \poids\vol^{\gmod} \geq
\int_{t\geq T} \left (\kappa + h^2 -\delta^2 + O(e^{-2t}) \right ) |\sigma|^2\poids\vol^{\gmod} ;
\end{equation}
compte tenu de l'inégalité $h\leq h_0 <0 $, le terme
$\kappa + h^2 -\delta^2$ est positif si
$$0\leq \delta <\sqrt {\kappa + h _0^2} ;
$$
  en choisissant $T$ suffisamment grand, on a  $\sqrt {\kappa +h_0^2}
>2$. 

Les métrique $\gmod$ et $\bar g$ diffèrent par un
$O(e^{-2t})$. En utilisant l'inégalité (\ref{isom00}) et le lemme
\ref{calculatoire} on en déduit ainsi le lemme
suivant~:
\begin{lemme}
\label{isom0}
Il existe  $\epsilon,c, T >0$  tels que pour tout
couple $(a,\phi)$ à support dans $\{t>T\}\times Y$ nul en $t=T$ et  $\delta \in
[0,2+\epsilon]$ on ait 
$$\int\langle P(a,\phi),(a,\phi)\rangle\poids\vol^{\bar g} \geq c\int \left (
|\nabla a|^2+ |\nabla^{B}\phi|^2 + |a|^2
+|\phi|^2 \right )\poids\vol^{\bar g},
$$
les normes étant prises par rapport à la métrique $\bar g$.
\end{lemme}

On peut maintenant facilement démontrer la proposition \ref{isomsw}.
\begin{demode}{de la proposition \ref{isomsw}}
D'après le lemme \ref{isom0} et la formule (\ref{exprPD}), on a pour
tout $\lambda
>0$ 
$$ \int \left ((1+\delta \lambda)|\mathcal
D(a,\phi)|^2 + \delta\lambda^{-1} ( |\phi|^2 + |a|^2)\right ) \vol \geq c \int
(|\nabla a|^2+ |\nabla^A\phi|^2+ |a|^2+|\phi|^2)\poids\vol; 
$$
en choisissant $\lambda$ suffisamment grand, on en déduit la proposition.
\end{demode}\medskip

\section{Démonstration du théorème~\ref{theo1}}
\label{secdemo}
\bigskip

On commence par une proposition similaire à celle du cas compact et qui est
la clef de la démonstration.
\begin{prop}
\label{plebrun}
Soit $M^4$ une variété orientée munie d'une métrique $g$
d'Einstein asymptotiquement hyperbolique  et d'une structure
$\spinc$ adaptée. Supposons que les équations
de Seiberg-Witten associées  admettent une
solution $(A,\Phi)$ telle que dans le théorème \ref{theo1}.
\begin{equation} 
\label{lebrun}
\int_M {s^2}\vol +  8 F_A\wedge F_A  \geq 0,
\end{equation}
avec égalité si et seulement si la métrique est de
Kähler-Einstein. Dans ce cas, on a  $\nabla_A\Phi=0$ et la
structure complexe se  déduit de $g$ et $F_A^+$.
\end{prop}\medskip
\begin{demo}
On procède comme sur les variétés compactes. D'après la
formule de Lichnerowicz et les équations de Seiberg-Witten,
$$ 0 = 4\langle \nabla^*_A\nabla_A \Phi,\Phi\rangle + s |\Phi|^2 +
|\Phi|^4.
$$
On multiplie cette identité par la fonction de troncature $\psi_{t_0}=
1-\chi_{t_0}$, où $\chi_{t_0}$ a été définie à la fin de la section
\ref{subhc}.
En appliquant l'inégalité de Cauchy-Schwarz, il vient
$$ \int \left (\langle \nabla_A\Phi, d\psi_{t_0}\otimes\Phi\rangle
+|\Phi |^4 \psi_{t_0} \right ) \vol \leq \left (\int s^2
\psi_{t_0} \vol \int |\Phi|^4 \psi_{t_0} \vol\right )^{\frac 12}  .
$$
Puisque $\nabla_A\Phi = O(e^{-(4+\epsilon)t})$  ce terme est $L^1$ et
on en déduit en utilisant la borne $C^0$ sur $\Phi$ que $\int \langle \nabla_A\Phi,
d\psi_{t_0}\otimes\Phi\rangle$ tend vers $0$ lorsque $t_0$ tend vers
l'infini. Il en résulte que
$$\int s^2 - |\Phi|^4 \vol \geq  0,
$$
avec égalité si et seulement si $\nabla_A\Phi=0$. D'après les
équations de Seiberg-Witten $|\Phi|^4 = 8|F_A^+|^2$. 
Puisque $F_A$ est imaginaire pure, $|F_A^+|^2\vol =  -F_A\wedge F_A +
|F_A^-|^2\vol$, et on en 
déduit l'inégalité de la proposition. Le cas d'égalité
entraîne  que $F_A^-=0$ et $\nabla^A\Phi=0$ et on obtient alors une
structure complexe à l'aide de la métrique et de la $2$-forme
autoduale parallèle $\sqrt 2 F_A/|F_A|$ pour laquelle la métrique
est kählérienne.
\end{demo}

Réciproquement, on a le lemme suivant~:
\begin{lemme}
Si $M$ admet une métrique de Kähler-Einstein  asymptotiquement
hyperbolique $\bar g$, alors les équations de Seiberg-Witten pour toute
métrique d'Einstein asymptotiquement symétrique $g$ (relativement à la même
structure de contact)  et la structure $\spinc$ canonique $\mathfrak s_0$ de
$\overline g$ admettent une solution.
\end{lemme}
\begin{demo}
La métrique de Kähler--Einstein définit un remplissage symplectique de
$M$ avec son bord de contact $Y$  à l'infini. D'après \cite{KM},
on a alors $SW(\mathfrak s_0)=1$ ; on en déduit une suite de solutions
$(A_\tau,\Phi_ \tau)$ des équations perturbées pour les métriques
$g_\tau$ approximant $g$. En appliquant le théorème~\ref{theoconv}, on obtient la
solution voulue.
\end{demo}

Finalement, sous les hypothèses du théorème $1$, on obtient pour $g$
une solution $(A,\Phi)$ des équation de Seiberg--Witten pour la structure
$\spinc$ canonique de $\CHH$. 
Si $B$ est la connexion de Chern de $K_{\CHH}^{-1}$ pour métrique hyperbolique
complexe, on sait alors que $F_{B}= \frac {is}4 \omega^{\CHH}$ et
$W^-_{\ghyp}=0$ d'où
\begin{equation}
\label{pl2}
\int \left(\shyp^2+48 |W^-_{\ghyp}|^2\right )\vol^{\mathcal H} + 8F_{B} \wedge F_{B}   = 0
\end{equation}
En choisissant une jauge où la décroissance est du type  $A-B=
O(e^{-(4+\epsilon)t})$, on voit en appliquant la formule de Stockes
que
$$ \int F_A\wedge F_A -  F_{B}\wedge F_{B}= 2\int d((A-B)\wedge F_{B})=0.
$$

 Pour des métriques
d'Einstein asymptotiquement hyperboliques, $ 2|W^-|^2 + \frac
{s^2}{24}$ est l'intégrant de $2\chi - 3\tau(M)$ et diverge mais  la différence entre les termes de bord pour la métrique $g$ et pour
$\ghyp$ tend vers $0$ d'après \cite{BH}. On en déduit d'après
(\ref{pl2}) que  
$$\int \left(s_g^2+48 |W^-_{g}|^2\right )\vol^g + 8F_{A} \wedge F_{A}
= 0,$$ 
ce qui entraîne d'après le lemme \ref{plebrun} que  la métrique $g$
est en réalité de Kähler--Einstein avec $W_g^-=0$.
Puisque $g$ est de Kähler--Einstein et autoduale il s'agit par
conséquent de la métrique hyperbolique complexe.
\section{Monopôles et structures CR}\bigskip

\label{ar}
Le théorème \ref{theoexistence} permet de montrer qu'il existe des
liens importants entre la géométrie d'une variété $Y^3$ munie d'une
structure CR et
les variétés $M^4$ munies de métriques d'Einstein asymptotiquement
hyperboliques complexes d'infini conforme $Y$. On dira que la
variété d'Einstein $M$ est un \emph{remplissage} de  $Y$.  
On a alors le corollaire suivant~:

\begin{cor}[Biquard, \cite{B3}] 
Soit $Y^3$ une variété munie d'une structure $CR$ et
 $M^4$ une variété munie d'une métrique  d'Einstein $g$ asymptotiquement
 hyperbolique d'infini conforme $Y$. Supposons de plus $M$ munie d'une
 structure $\spinc$ adaptée $\mathfrak s$ avec un invariant de
 Kronheimer--Mrowka $SW(\mathfrak s)\neq 0$. 
On a alors l'inégalité de Miyaoka--Yau
\begin{equation}
\label{rbh}
0\leq \chi(M)-3\tau(M)+\nu(\partial_\infty M) = \frac 1{8\pi^2} \int_M
\left ( 3|W_g^-|^2- |W_g^+|^2+ 
\frac {s_g^2}{24}\right ) \vol^g   , 
\end{equation}
où $\nu$ est l'invariant de la structure CR  à l'infini défini
dans \cite{BH}.
De plus on a égalité  si et seulement si la métrique $g$ est hyperbolique
complexe.
\end{cor}

Notons que l'égalité de (\ref{rbh}) est établie dans \cite{BH}.

\begin{demo}
Par hypothèse l'invariant de Seiberg--Witten est non nul ; par conséquent la
dimension virtuelle de l'espace des modules est positive ou nulle soit
(cf. \cite{KM}, th. $2.4$)
$$d(\mathfrak s) = \langle e(W^+,\Phi_0),[X,\partial X]\rangle \geq 0,
$$
avec $ e(W^+,\Phi_0)$ la classe d'Euler relative du fibré des spineurs
$W ^+$. 

On calcule maintenant ce nombre caractéristique via la théorie
de Chern--Weil. Soit $A_0$ la connexion de Chern du fibré
anti-canonique induite par  la métrique de
Kähler--Einstein $\bar g$ définie près de l'infini. On étend $A_0$ et
$\bar g$
arbitrairement au dessus de $M$ tout entier par partition de l'unité. 
On en déduit une connexion
spinorielle $\nabla_{A_0}$ sur $W^+$ telle que $\nabla_{A_0}\Phi_0 =
0$ près de l'infini. En ce qui concerne la construction de la
classe de Thom $U$ de $W^+$, on  se référera à  \cite{BGV},
section 1.6.  Soit $\Psi$ un spineur de $W^+$ transverse à la section
nulle et égal à $\Phi_0$ près de
l'infini. On a la formule de transgression
\begin{equation}
\label{trans}
 \mathrm{Pf}(R(\nabla_{A_0})) - \Psi^* U  =  -i d\int_0^1
T \left (\Psi\wedge e^{-(r^2|\Psi|^2/2) + ir\nabla_{A_0}\Psi +
R(\nabla_{A_0})} \right ) dr,
\end{equation}
où $T : \Omega (M, \Lambda W^+ ) \rightarrow \Omega (M) $ est
l'intégrale de Berezin et la courbure $R(\nabla_{A_0})$ a été
identifiée à un élément de $\Omega^2(M,\Lambda^2W^+)$. Comme $\Phi_0$
est parallèle, on en déduit que l'intégrale de (\ref{trans}) est à
support compact. Par conséquent,
\begin{align*}
    \langle e(W^+,\Phi_0),[X,\partial X]\rangle & =  \int_M
\mathrm{Pf}(R(\nabla_{A_ 0})) \vol^{\bar g } \\
&= \frac 1{4\pi^2} \int_M
-F_{A_0}\wedge F_{A_0} - \left ( 2|W_{\bar g}^+|^2 - \frac 12
|\Ric^{\bar g}_0|^2+ \frac
{s_{\bar g}^2}{24}\right ) \vol^{\bar g}.
\end{align*}

Soit $(A,\Phi)$ une solution des équations de
Seiberg--Witten fournie par le théorème \ref{theoexistence}. On en
déduit une connexion spinorielle $\nabla_A$ grâce à la métrique $g$ ;
alors en utilisant la décroissance de $A-A_0$ et  la formule de Stockes
on obtient $\int F_{A}\wedge F_{A} - F_{A_0}\wedge F_{A_0} =0$. Par
ailleurs, d'après \cite{BH}, on a
 $$\int_M
 \left ( 2|W_{ g}^+|^2+ \frac
{s_{ g}^2}{24}\right ) \vol^{ g}
 - \left ( 2|W_{\bar g}^+|^2 - \frac 12
|\Ric^{\bar g}_0|^2+ \frac
{s_{\bar g}^2}{24}\right ) \vol^{\bar g} = 0.$$
Il en résulte que 
$$\frac 1{4\pi^2} \int_M -F_A\wedge F_A - \left ( 2|W_g^+|^2 + \frac
{s_g^2}{24}\right ) \vol^g = d(\mathfrak s) \geq 0.
$$
En utilisant la proposition \ref{plebrun}, on obtient l'inégalité
$$0\leq \frac 1{4\pi^2}\int_M \left ( \frac {s_g^2}{8} -2 |W_g^+|^2- \frac
{s_g^2}{24} \right )\vol ^g \leq \frac 1{2\pi^2}\int_M \left ( 3|W_g^-|^2- |W_g^+|^2+ \frac
{s_g^2}{24}\right )\vol ^g,
$$
avec égalité si et seulement si $W_g^-=0$ et la métrique $g$ est
Kähler--Einstein, ce qui implique que la métrique $g$ est en fait
hyperbolique complexe. 
\end{demo}


\bigskip


\begin{thebibliography}{MMM}


 \bibitem[Bi]{B}  O. Biquard,  \emph{Métriques d'Einstein asymptotiquement
symétriques}, Astérisque {\bf 265}  (2000).

 \bibitem[Bi2]{B2} O. Biquard, \emph{Métriques d'Einstein à cusps et
     équations de Seiberg-Witten,} J. Reine Angew. Math. {\bf 490},
     129--154  (1997).

 \bibitem[Bi3]{B3} O. Biquard, communication privée (déc. 2001).


 \bibitem[BiH]{BH}  O. Biquard, M. Herzlich, \emph{A Burns-Epstein
 invariant for ACHE manifolds}, arXiv math.DG/0111218 (2001).

 \bibitem[BCG]{BCG}  G. Besson, G. Courtois, S. Gallot,
     \emph{Entropies et rigidités des espaces localement symétriques
     de courbure strictement négative,} Geom. Func. Anal.
 {\bf 5}, 731--799  (1995). 


 \bibitem[BG]{BG}  J.P. Bourguignon, P. Gauduchon, \emph{Spineurs,
     opérateurs de Dirac et variations de métriques,} Comm. Math. Phys.
 {\bf 144}, 581--599  (1992). 

 \bibitem[BGV]{BGV}  N. Berline, E. Getzler, M. Vergne,  \emph{Heat
     Kernels and Dirac Operators,} Grund. Math. 
 {\bf 298}, Springer--Verlag. 


 \bibitem[BoH]{BoH}  H. Boualem, M. Herzlich, \emph{Rigidity at infinity for even-dimensional asymptotically complex hyperbolic spaces}, Prep. Univ. Montpellier II
 {\bf 10} (2001).



 \bibitem[CY]{CY}  S.Y. Cheng, S.T. Yau, \emph{On the existence of a
 complete Kähler metric on non--compact complex manifolds and the
 regularity of Fefferman's equation,}  Comm. Pure. App. Math {\bf 33},
 507--544  (1980).  








 \bibitem[H]{H}   M. Herzlich, \emph{Scalar curvature and rigidity for
 odd-dimensional complex hyperbolic spaces}, Math. Ann. {\bf 312}, 641--657 (2001).


 \bibitem[KM]{KM}  P.B. Kronheimer, T.S. Mrowka, \emph{Monopoles and
  contact structures,} Invent. Math. { \bf 130},  209--255 (1997).
 

 \bibitem[KM2]{KM2}  P.B. Kronheimer, T.S. Mrowka, \emph{The genus of
  embedded surfaces in the projective plane,} Math. Res. Lett.  { \bf 1},  797--808 (1994).


 \bibitem[Ko]{Ko}  D. Kotschick, \emph{The Seiberg-Witten invariants
  of symplectic four-manifolds, [after C.H. Taubes],} Sém. N. Bourbaki { \bf 812},  (1995-1996).

 \bibitem[L]{L}  C. Le\,Brun, \emph{Einstein Metrics and Mostow
Rigidity,} Math. Res. Lett. { \bf 2}, 1--8 (1995).


 \bibitem[R]{R2}  Y. Rollin, \emph{Rigidité d'Einstein du plan
 hyperbolique complexe}, C. R. Acad. Sci. Paris, Ser. I {\bf 334},
 671--676 (2002). 


 \bibitem[R2]{R}  Y. Rollin, \emph{Surfaces kählériennes de volume fini
 et équations de Seiberg--Witten}, arXiv math.DG/0106077 (2001)
 ---  à paraître dans Bull. Soc. Math. Fr.

 \end{thebibliography}
\end{document}